\documentclass[11pt]{article}   
\begin{document}
\newtheorem{defn0}{Definition}[section]
\newtheorem{prop0}[defn0]{Proposition}
\newtheorem{thm0}[defn0]{Theorem}
\newtheorem{lemma0}[defn0]{Lemma}
\newtheorem{coro0}[defn0]{Corollary}
\newtheorem{exa}[defn0]{Example}
\newtheorem{exe}[defn0]{Exercise}
\newtheorem{notation}[defn0]{Notation}
\newtheorem{remark}[defn0]{Remark}
\newtheorem{question}[defn0]{Question}
\def\rig#1{\smash{ \mathop{\longrightarrow}\limits^{#1}}}
\def\swar#1{\swarrow\rlap{$\vcenter{\hbox{$\scriptstyle#1$}}$}}
\def\lswar#1{\swarrow\llap{$\vcenter{\hbox{$\scriptstyle#1$}}$}}
\def\sear#1{\searrow\rlap{$\vcenter{\hbox{$\scriptstyle#1$}}$}}
\def\lsear#1{\searrow\llap{$\vcenter{\hbox{$\scriptstyle#1$}}$}}   
\def\near#1{\nearrow\rlap{$\vcenter{\hbox{$\scriptstyle#1$}}$}}   
\def\dow#1{\Big\downarrow\rlap{$\vcenter{\hbox{$\scriptstyle#1$}}$}}
\def\ldow#1{\Big\downarrow\llap{$\vcenter{\hbox{$\scriptstyle#1$}}$}}   
\def\up#1{\Big\uparrow\rlap{$\vcenter{\hbox{$\scriptstyle#1$}}$}}
\def\lef#1{\smash{ \mathop{\longleftarrow}\limits^{#1}}}
\newcommand{\defref}[1]{Definition~\ref{#1}}
\newcommand{\propref}[1]{Proposition~\ref{#1}}
\newcommand{\thmref}[1]{Theorem~\ref{#1}}
\newcommand{\lemref}[1]{Lemma~\ref{#1}}
\newcommand{\corref}[1]{Corollary~\ref{#1}}
\newcommand{\exref}[1]{Example~\ref{#1}}
\newcommand{\secref}[1]{Section~\ref{#1}}
\newcommand{\qedd}{\hfill\framebox[2mm]{\ }}
\def\P#1{{\bf P}^#1}
\def\I{{\cal I}}
\def\O{{\cal O}}
\def\C{{\bf C}}
\def\Z{{\bf Z}}
\def\proof{{\it Proof\ }}
\def\M{{\cal M}}
\def\F{{\cal F}}
\def\Q{{\bf Q}}
\def\N{{\bf N}}
\long\def\ignore#1{}
\title{Induction for secant varieties of Segre varieties}

\date{May 26th, 2006}

\author{Hirotachi Abo
 \and Giorgio Ottaviani
 \and Chris
Peterson}

\date{\ }

\maketitle

\vspace*{-0.75in}
\begin{abstract}

\noindent This paper studies the dimension of secant varieties to Segre varieties. The problem is cast both in the setting of
tensor algebra and in the setting of algebraic geometry. An inductive procedure is built around the ideas of successive specializations of points and projections. This reduces the calculation of the dimension of the secant variety in a high dimensional case
to a sequence of calculations of partial secant varieties in low dimensional
cases. As applications of the technique: We give a complete classification of defective $t$-secant varieties to Segre varieties for $t\leq 6$. We generalize a theorem of Catalisano-Geramita-Gimigliano on non-defectivity of tensor powers of $\P {n}$. We determine the set of $p$ for which unbalanced Segre varieties have defective $p$-secant varieties. In addition, we completely describe the dimensions of the secant varieties to the deficient Segre varieties $\P {1}\times \P {1} \times \P {n} \times \P {n}$ and $\P {2}\times \P {3} \times \P {3}$. In the final section we propose a series of conjectures about defective Segre varieties.

\vskip 3pt


\vskip 3pt

\end{abstract}

\tableofcontents

\

\footnotetext[1]{The authors would like to thank the Department of Mathematics of the Universit\`a di Firenze, the Department of Mathematics of Colorado State University, the National Science Foundation and GNSAGA of Italian INDAM.}

\pagebreak
\section{Introduction}
If $Q_1,\dots, Q_p$ are points then we let $<Q_1,\dots, Q_p>$ denote their linear span. Let $X_1,\dots, X_p\subseteq\P m$ be projective varieties of dimensions $d_1,\dots, d_p$. The {\it join} of the varieties, $J(X_1,\dots, X_p)$, is defined to be the Zariski closure of the union of the linear span of $p$-tuples of points $(Q_1,\dots, Q_p)$ where $Q_i\in X_i$.  In other words $$J(X_1,\dots, X_p)=\overline{\bigcup_{Q_1\in X_1,\dots, Q_p\in X_p}<Q_1,\dots,Q_p>}.$$ The expected
dimension (and the maximum possible dimension) of $J(X_1,\dots, X_p)$ is $\min\{m,p-1+\sum d_i\}$. If $X\subseteq \P m$ is a variety then the $p$-secant variety of $X$ is defined to be the join of $p$ copies of $X$. We will denote this by $\sigma_p(X)$. Hence $\sigma_1(X)=J(X)=X$ while $\sigma_2(X)=J(X,X)$ is the variety of secant lines to $X$. The expected
dimension (and the maximum possible dimension) of $\sigma_p(X)$ is $\min\{m,pr+(p-1)\}$. $X$ is said to have a {\it defective $p$-secant variety} if $\dim\sigma_p(X)<\min\{m,pr+(p-1)\}$. $X$ is called {\it defective} if there exists a $p$ such that
$\dim\sigma_p(X)<\min\{m,pr+(p-1)\}$. In other words, $X$ is defective if for \underline{some} $p$, $X$ has a defective $p$-secant variety. For instance, a classical theorem in algebraic geometry states that the Veronese surface $V\subset \P{5}$ is defective since the dimension of $\sigma_2(V)$ is 4 (instead of the expected dimension of 5).

Let $\P{{n_i}}={\bf P}(V_i)$ where $V_i$ is a vector space of dimension $n_i+1$ over a field of characteristic zero, not necessarily algebraically closed.
The aim of this note is to compute the dimension of $\sigma_p(X)$ when
$X$ is a Segre variety $\P {{n_1}}\times\ldots\times\P{{n_k}}$ embedded in
${\bf P} (V_1\otimes\ldots\otimes V_k)$. 
We say that $(n_1,\ldots ,n_k)$ is defective if there exists a $p$ such that
$\dim\sigma_p(\P {{n_1}}\times\ldots\times\P{{n_k}})$ is less than the expected dimension $\min\left\{\prod(n_i+1)-1, s(\sum n_i)+s-1 \right\}$. If $W_1,\dots, W_p \subseteq X\subseteq \P m$, then $J(W_1,\dots ,W_p)$ is called a partial secant variety of $X$.
In Section 2, we describe the basic tensor algebra that will be used throughout the paper. In Section 3, we give an inductive procedure that reduces the computation of $\dim\sigma_p(\P {{n_1}}\times\ldots\times\P{{n_k}})$ to the computation of the dimension of a collection of partial secant varieties of low dimensional Segre varieties. Thus, a high dimensional problem is reduced, inductively, to a collection of easily computable low dimensional problems. In Section 4, we apply this procedure to give a complete classification of defective $t$-secant varieties to Segre varieties for $t\leq 6$.  In the process of carrying out the classification, we characterize the set of $p$ for which unbalanced Segre varieties have defective $p$-secant varieties. Modulo the unbalanced Segre varieties, there seem to be very few defective cases. However, we show that the Segre varieties $\P {1}\times \P {1} \times \P {n} \times \P {n}$ and $\P {2}\times \P {3} \times \P {3}$ are defective (and completely describe the dimensions of their secant varieties). In Section 5, we generalize a theorem of Catalisano-Geramita-Gimigliano on the non-defectivity of tensor powers of $\P {n}$. We close the paper with a series of conjectures on the existence and classification of defective Segre varieties. In addition to evidence provided by the theorems of this paper, further evidence in support of the conjectures can be obtained via Montecarlo techniques in a computer algebra system such as CoCoA, Macaulay 2 or  Singular \cite {Co,GS,GPS05}.

The interest in this subject comes from several different sources. In algebraic geometry, the Segre varieties form an important class of geometric objects. In one guise, points on a Segre variety, $V$,  are viewed as parametrizing rank one (or decomposable) tensors. A tensor is said to have {\it rank $r$} if it can be written as a linear combination of $r$ rank one tensors (but not fewer). A tensor is said to have {\it border rank $r$} if it can be expressed as the limit of rank $r$ tensors but not as the limit of rank $r-1$ tensors. With this notation, $\sigma_p(V)$ parametrizes tensors with border rank at most $p$. Alternatively, these same ideas can be expressed in terms of decomposition of multidimensional matrices as linear combinations of simpler ``rank 1" multidimensional matrices (\cite{GKZ}, 
\cite{CGG1}). In numerical analysis  a thorough understanding of the dimension of $\sigma_p(V)$ has applications to complexity theory, for example to algorithms for matrix 
multiplication (\cite{BCS}, \cite{La}). More recently this topic appears, through its relationship with algebraic statistics and higher order correlations, in connection with 
computational biology (\cite{ERSS}).
The special case $X=\P 1\times\ldots\times\P 1$ (\cite{CGG2}) has made several appearances in the recent physics
literature (see for example \cite{LT} and the literature quoted therein). The interested reader should also consider the accessible articles (\cite{BM},\cite{C}) for an overview of some related topics.

\section{Basic tensor algebra for Segre Varieties}

In this section, questions about secant varieties to Segre varieties are reinterpreted as questions in tensor algebra. We begin by introducing the notation that will be used throughout this paper.

\begin{defn0} Let $Y$ be a subspace of a vector space $V$. Let $V^{\vee}$ denote the dual vector space of $V$. The orthogonal, $Y^{\perp}$ of $Y$, is defined by
$$Y^{\perp}:=\{\omega\in V^{\vee}\hskip 1pt | \hskip 1pt \omega(v)=0\quad\forall\ v\in Y\}.$$
\end{defn0}
It is worth noting that the dimension of $Y$ in $V$ is the same as the codimension of $Y^{\perp}$ in $V^{\vee}$.

The symmetric algebra of a vector space $V$, $Sym(V)=\bigoplus_{i=0}^{\infty}Sym^i(V)$, comes equipped with a natural grading.
Let $S(V_i)={\bf C}\oplus V_i$ be the truncated symmetric algebra arising as
the quotient of the symmetric algebra by the ideal of elements whose degree is greater than or equal to 2 (in the natural grading).
Given vector spaces $V_1, \dots ,V_k$, the commutative algebra $T=S(V_1)\otimes\ldots\otimes S(V_k)$ has a multi-gradation
indexed by $k$-tuples of non-negative integers where the summand corresponding
to $n=(n_1,\ldots ,n_k)$ is zero if some $n_i\ge 2$. We will let $T_{n_1,\ldots ,n_k}$ denote the summand of $T$ with multi-degree $(n_1,\ldots ,n_k)$. In particular  $T_{0,\ldots, 1,\ldots, 0}=V_i$ and
$T_{1,\ldots, 1}=V_1\otimes\ldots \otimes V_k$ are direct summands of $T$ with multi-degrees $(0,\ldots, 1,\ldots, 0)$ and $(1,\ldots, 1)$ respectively. 
Since $(A\otimes B)^{\vee}=A^{\vee}\otimes B^{\vee}$, we have $$T^{\vee}=S(V_1^{\vee})\otimes\ldots\otimes 
S(V_k^{\vee}),\hskip 10pt T^{\vee}_{0,\ldots, 1,\ldots, 0}=V_i^{\vee}
 \hskip 20pt {\rm and}\hskip 20pt T^{\vee}_{1,\ldots ,1}=V_1^{\vee}\otimes\ldots \otimes V_k^{\vee}.$$

Let $<v_i>^{\perp}$ denote the homogeneous ideal in $T^{\vee}$ which is generated
by the subspace $<v_i>^{\perp}\hskip 2pt\subseteq V_i^{\vee}$. Though $<v_i>^{\perp}$ denotes both a homogeneous ideal and a subspace, in this paper there will be no danger of ambiguity. The following lemma is analogous to the well known cases of
projective spaces and Grassmann varieties \cite{CGG3}.

\begin{lemma0}
\label{tangent} Let $p=v_1\otimes\ldots\otimes v_k$ be a point of
$X={\bf P}(V_1)\times\ldots\times{\bf P}(V_k)$. Then 

$(i)\qquad T_pX=V_1\otimes v_2\otimes\ldots\otimes v_k+v_1\otimes V_2\otimes\ldots\otimes v_k+\ldots
 +v_1\otimes v_2\otimes\ldots\otimes V_k$

$(ii)\qquad T_pX^{\perp}= \left[(<v_1>^{\perp}+\ldots+<v_k>^{\perp})^2
\right]_{1,\ldots, 1}\subseteq V_1^{\vee}\otimes\ldots\otimes V_k^{\vee}$.
\end{lemma0}

\proof\hskip -3pt . To see (i),  take the derivative of the parametric curve
$(v_1+\epsilon v_1')\otimes\ldots\otimes (v_k+\epsilon v_k')$ at $\epsilon=0$ and let $v_1', v_2',\dots, v_k'$ vary over $V_1,V_2, \dots, V_k$.

In order to prove (ii), consider that
$$\left(v_1\otimes\ldots \otimes v_{i-1}\otimes V_i\otimes v_{i+1}\otimes \ldots\otimes v_k\right)^{\perp}=\left(\sum_{j\neq i}<v_j>^{\perp}\right)_{1,\ldots, 1}, \hskip 10pt {\rm hence}$$
$$T_pX^{\perp}=\bigcap_{i=1}^k\left(\sum_{j\neq i}<v_j>^{\perp}\right)_{1,\ldots, 1}.$$

Complete $v_i=v_{i,1}$ to a basis $\{v_{i,1},\ldots, v_{i,n_i+1}\}$ of $V_i$.
We label the dual basis of $V_i^{\vee}$ by $\{v^{i,1},\ldots, v^{i,n_i+1}\}$. In the dual basis,
$(<v_j>^{\perp})$ is generated by $\{v^{j,2},\ldots, v^{j,n_j+1}\}$.

Now $\left(\sum_{j\neq i}<v_j>^{\perp}\right)_{1,\ldots, 1}$ contains all monomials with multidegree $(1,\dots, 1)$  with the exception of the following $n_i+1$
$$\{v^{1,1}\otimes v^{2,1}\otimes\ldots \otimes v^{i-1,1}\otimes v^{i,j}\otimes v^{i-1,1}\otimes\ldots
\otimes v^{k,1}\ |\ 1\leq j\leq n_i+1\}.$$

Hence $\bigcap_{i=1}^k\left(\sum_{j\neq i}<v_j>^{\perp}\right)_{1,\ldots, 1}$ is generated by all basis elements
$\alpha_1\otimes\ldots\otimes\alpha_k$ with $\alpha_j\neq v^{j,1}$ for at least two different values of the index $j$. These are
exactly the generators of $\left[(<v_1>^{\perp}+\ldots+<v_k>^{\perp})^2
\right]_{1,\ldots, 1}$.
\qedd

\

A subspace $Y\subseteq V_1\otimes\ldots \otimes V_k$ is called {\it monomial} if
there exist bases of $V_1,\ldots ,V_k$ such that a basis of $Y$ can be
expressed in terms of monomials in the bases of $V_1,\ldots ,V_k$.

\begin{coro0} Let $p=v_1\otimes\ldots\otimes v_k$ be a point of
$X={\bf P}(V_1)\times\ldots\times{\bf P}(V_k)$. Then $T_pX$ and $T_pX^{\perp}$ are monomial subspaces.
\end{coro0}

Fix now a  subspace $H\subseteq V_1$  of dimension $h$.
For any $p=v_1\otimes\ldots\otimes v_k$,  we have either
$v_1\notin H$ or $v_1\in H$.

\medskip
Let $$f\colon V_1^{\vee}\otimes V_2^{\vee}\otimes\ldots\otimes V_k^{\vee}\rig{}
H^{\vee}\otimes V_2^{\vee}\otimes\ldots\otimes V_k^{\vee}$$ be the natural projection and let 
\begin{equation}\label{induction}0\rig{}K\rig {}T_pX^{\perp}\rig{} f\left(T_pX^{\perp}\right) \rig{}0
\end{equation}
be the restriction exact sequence, where
$K=T_pX^{\perp}\cap \left[(V_1/H)^{\vee}\otimes V_2^{\vee}\otimes\ldots\otimes V_k^{\vee}\right]$.
Clearly both $f\left(T_pX^{\perp}\right)$ and $K$ depend heavily on whether $v_1\notin H$ or $v_1\in H$. This dependence is captured in the following:

\begin{lemma0}
\label{restriction} Consider a point $v_1\in V$ and a subspace $H\subseteq V$.
\begin{enumerate}
        \item  
If  $v_1\notin H$, let $[v_1]\in V/H$ denote its quotient class. We have $$f\left(T_pX^{\perp}\right)= \left[<v_2>^{\perp}+\ldots+<v_k>^{\perp}\right]_{1,\ldots, 1}$$
which has codimension $h$ in $H^{\vee}\otimes V_2^{\vee}\otimes\ldots\otimes V_k^{\vee}$,
and  
$$K=\left[(<[v_1]>^{\perp}+<v_2>^{\perp}+\ldots+<v_k>^{\perp})^2\right]_{1,\ldots, 1}$$
which has codimension $1+\sum_{i=2}^kn_i+(n_1-h)$ in $(V_1/H)^{\vee}\otimes V_2^{\vee}\otimes\ldots\otimes V_k^{\vee}$.
 
 \item
If  $v_1\in H$, we have $$f\left(T_pX^{\perp}\right)= \left[(<v_1>^{\perp}+<v_2>^{\perp}+
\ldots+<v_k>^{\perp})^2\right]_{1,\ldots, 1}$$
which has codimension $h+\sum_{i=2}^kn_i$ in $H^{\vee}\otimes V_2^{\vee}\otimes\ldots\otimes V_k^{\vee}$,
and 
$$K=\left[<v_2>^{\perp}+\ldots+<v_k>^{\perp}\right]_{1,\ldots, 1}$$
which has codimension $n_1+1-h$ in $(V_1/H)^{\vee}\otimes V_2^{\vee}\otimes\ldots\otimes V_k^{\vee}$.
\end{enumerate}
\end{lemma0}

\proof\hskip -3pt .  We first consider the case where $v_1 \notin H$. In this setting,
$(<v_1>^{\perp})$ projects to the entire subspace $H^{\vee}$. Hence every element 
in $\left[<v_2>^{\perp}+\ldots+<v_k>^{\perp}\right]_{1,\ldots, 1}$ is the projection of
an element of $(<v_1>^{\perp})\hskip 2pt\cap\hskip 2pt (<v_i>^{\perp})$ for some $i$.
Both the assertion about $K$ and the inclusion $f\left(T_pX^{\perp}\right)\subseteq  \left[<v_2>^{\perp}+\ldots+<v_k>^{\perp}\right]_{1,\ldots, 1}$ are clear. From (ii) of \lemref{tangent}, we have $\left[<v_2>^{\perp}+\ldots+<v_k>^{\perp}\right]_{1,\ldots, 1}\subseteq f\left(T_pX^{\perp}\right) $.
The proof for the case where $v_1\in H$ is analogous and is left to the reader. Note that $T_pX^{\perp}$ has codimension $1+\sum_{i=1}^kn_i$ in 
$V_1^{\vee}\otimes\ldots \otimes V_k^{\vee}$. From this fact, the statements about the codimension of $K$ follow.
\qedd

\

We now look at the connection with the secant varieties of $X=\P {{n_1}}\times\ldots\times\P{{n_k}}$.
The expected dimension of $\sigma_p(X)$ is 
$$\min\{ \prod_{i=1}^k(n_i+1)-1,p(\sum_{i=1}^kn_i)+(p-1)\}.$$

 There is a unique integer $s$ such that
$\sigma_s(X)$ fills the ambient space 
and $\sigma_{s-1}(X)$ does not. The expected value for such an $s$ is
$$S( n_1,\ldots, n_k ):= \left\lceil \frac{\prod_{i=1}^k(n_i+1)}
{(\sum_{i=1}^kn_i)+1}\right\rceil.$$

A standard application of Terracini's lemma, as in \cite{CGG2}, shows that
$\sigma_s(X)$ has the expected dimension
 if and only if for
$s$ generic points $p_1,\dots, p_s$,  the linear space $\left[ T_{p_1}X^{\perp}\cap\ldots\cap
T_{p_s}X^{\perp}\right]$ has the expected codimension in
 $T^{\vee}_{1,\ldots, 1}$, that is

$$
\left\{\begin{array}{lcc} {s(\sum_{i=1}^kn_i)+1}\
&\hbox{\ for\ }&s<S(n_1,\ldots, n_k )\\
\prod_{i=1}^k(n_i+1)&\hbox{\ for\ }&s\ge S(n_1,\ldots, n_k )\\
\end{array}\right.$$ 

Consider again the point $p=v_1\otimes\ldots\otimes v_k$.
\lemref{restriction} suggests we focus our attention on the subspaces
 $G^i_pX\subseteq T_{1,\ldots ,1}$ defined by
 
 $$G^i_pX^{\perp}=\left[(\sum_{j\neq i}<v_j>^{\perp})
\right]_{1,\ldots, 1}.$$
It is easy to check that
$$G^i_pX=\left(v_1\otimes\ldots \otimes v_{i-1}\otimes V_i\otimes v_{i+1}\otimes\ldots\otimes v_k\right)
 $$
has dimension $n_i+1$ in $T_{1,\ldots, 1}$ (and that
$G^i_pX^{\perp}$ has codimension $n_i+1$ in $T_{1,\ldots, 1}^{\vee}$). 

\begin{remark} We sketch the geometrical construction which is behind the tensor algebra of this section. We have denoted by ${\bf P}(V_1)$ the projective space of 
lines
in $V_1$, so that $H^0({\bf P}(V_1),\O(1))=V_1^{\vee}$. The subvariety $X'={\bf P}(H)\times{\bf P}(V_2)\times
\ldots\times{\bf P}(V_k)\subset X$ is the zero locus of a section of the vector bundle
$(V_1/H)\otimes\O_X(1,0\ldots, 0)$. We get the Koszul complex
\begin{equation}\label{behind}
\ldots\rig{}\wedge^2(V_1/H)^{\vee}\otimes\O_X(-2,0\ldots, 0)
\rig{}(V_1/H)^{\vee}\otimes\O_X(-1,0\ldots, 0)\rig{}\O_X\rig{}\O_{X'}\rig{}0
\end{equation}
After tensoring {\rm (\ref{behind})} by   $\O_X(1,1\ldots, 1)$ and taking cohomology we get
$$0\rig{}  (V_1/H)^{\vee}\otimes\ V_2^{\vee}\otimes\ldots\otimes V_k^{\vee}
\rig{} V_1^{\vee}\otimes\ V_2^{\vee}\otimes\ldots\otimes V_k^{\vee}\rig{}
H^{\vee}\otimes\ V_2^{\vee}\otimes\ldots\otimes V_k^{\vee}\rig{}0$$
Let $p$ be a double point on X.
After tensoring {\rm (\ref{behind})} by   $I_{p}^2\otimes\O_X(1,1\ldots, 1)$ and taking cohomology we get
exactly sequence {\rm (\ref{induction})}:

\small{\[\begin{array}{ccccccccc}
0&\rightarrow&K&\rightarrow& T_pX^{\perp}&\rightarrow&f(T_pX^{\perp})&\rightarrow&0\\
&&\cap&&\cap&&\cap\\
0&\rightarrow&  (V_1/H)^{\vee}\otimes\ V_2^{\vee}\otimes\ldots\otimes V_k^{\vee}
&\rightarrow& V_1^{\vee}\otimes\ldots\otimes V_k^{\vee}&\rightarrow &
H^{\vee}\otimes\ V_2^{\vee}\otimes\ldots\otimes V_k^{\vee}&\rightarrow&0
\end{array}\]}

Hence, in the language of \cite{AH} 
$f(T_pX^{\perp})$ plays the role of  trace and $K$ plays the role of residual.
\end{remark}

\section{Induction for secant varieties to Segre varieties}
In this section, we develop a method of induction for secant varieties to Segre varieties. 

\begin{notation} We fix now the notation that will be used throughout this section.
\begin{itemize}
\item $X=\P {{n_1}}\times\P{{n_2}}\times\dots \times\P{{n_k}}$
\item If $\vec {\mathbf n}=(n_1,\dots, n_k)$ then ${\mathbf P}^{\vec {\mathbf n}}=\P {{n_1}}\times\P{{n_2}}\times\dots \times\P{{n_k}}$
\item For $s$ generic points $p_1,\ldots p_s\in X$ let $T_sX=T_{p_1}X+\ldots +T_{p_s}X$
\item For $t$ generic points $q_1,\ldots ,q_t\in X$
let $G^i_tX=G^i_{q_1}X+\ldots + G^i_{q_t}X$
\end{itemize}
\end{notation}

This notation leads to the following fundamental definition.

\begin{defn0} Let $s,a_1,a_2,\dots, a_k$ be non-negative integers and let $X={\mathbf P}^{\vec {\mathbf n}}$.
\begin{itemize}
\item If for $s+a_1+a_2+\dots +a_k$ generic points, the linear space spanned by
$T_sX+G^1_{a_1}X+G^2_{a_2}X+\dots +G^k_{a_k}X\subseteq T_{(1,\ldots ,1)}$  has dimension
$$D=\min\{s(1+\sum_{i=1}^k{n_i})+\sum_{i=1}^k (a_i(n_i+1)),\prod_{i=1}^k{(n_i+1)}\}$$ then we say that $T(n_1, \dots, n_k;s ; a_1,\dots , a_k)$ is true. At times we will abbreviate this by $T(\vec {\mathbf n},s,\vec {\mathbf a})$. By duality, we have the equivalent definition  that $T(\vec {\mathbf n},s,\vec {\mathbf a})$ is true if and only if for  $s+\sum a_i$
generic points, the intersection  $T_sX^{\perp}\cap G^1_{a_1}X^{\perp}\cap G^2_{a_2}X^{\perp}\cap\dots \cap G^k_{a_k}X^{\perp}
\subseteq T^{\vee}_{(1,\ldots ,1)}$
has  codimension D.
\item If $s(1+\sum{n_i})+\sum (a_i(n_i+1))\leq \prod{(n_i+1)}$ then $(\vec {\mathbf n},s,\vec {\mathbf a})$ is called {\bf subabundant}.
\item If $s(1+\sum{n_i})+\sum (a_i(n_i+1))\geq \prod{(n_i+1)}$ then $(\vec {\mathbf n},s,\vec {\mathbf a})$ is called {\bf superabundant}.
\item If $s(1+\sum{n_i})+\sum (a_i(n_i+1))= \prod{(n_i+1)}$ then $(\vec {\mathbf n},s,\vec {\mathbf a})$ is called {\bf equiabundant}.
\item If $(\vec {\mathbf n},s,\vec {\mathbf 0})$ is equiabundant and $T(\vec {\mathbf n},s,\vec {\mathbf 0})$ is true then ${\mathbf P}^{\vec {\mathbf n}}$ is called {\bf perfect}.
\item If $(\prod_{i=1}^k(n_i+1))/(1+\sum_{i=1}^kn_i)$ is an integer then $\vec {\mathbf n}$ is called {\bf numerically perfect}.
\end{itemize}
\end{defn0}

For efficiency, we will often write statements such as $T(\vec {\mathbf n},s,\vec {\mathbf a})$ is true and subabundant when we really mean $T(\vec {\mathbf n},s,\vec {\mathbf a})$ is true and $(\vec {\mathbf n},s,\vec {\mathbf a})$ is subabundant.

\begin{remark} \label{usefulrem} Given two $k$-dimensional vectors $\vec {\mathbf n}, \vec {\mathbf n}'$, we say $\vec {\mathbf n}'\leq \vec {\mathbf n}$ if $n_i'\leq n_i$ for each $1\leq i\leq k$. We make three simple remarks:
\begin{itemize}
\item[(i)]  $T(n_1, \dots, n_k;s ; 0, \dots, 0)$ is true if and only if 
$\sigma_s( \P {{n_1}}\times\P{{n_2}}\times\dots\times\P{{n_k}})$ has the expected dimension.
\item[(ii)] If $T(\vec {\mathbf n},s,\vec {\mathbf a})$ is true and subabundant then $T(\vec {\mathbf n},s',\vec {\mathbf a}')$ is true and subabundant for any choice of $s',\vec {\mathbf a}'$ with
$s'\leq s$ and  $\vec {\mathbf a}'\leq\vec {\mathbf a}$.
\item[(iii)] If $T(\vec {\mathbf n},s,\vec {\mathbf a})$ is true and superabundant then $T(\vec {\mathbf n},s',\vec {\mathbf a}')$ is true and superabundant for any choice of $s',\vec {\mathbf a}'$ with
$s\leq s'$ and  $\vec {\mathbf a}\leq\vec {\mathbf a}^\prime$.
\end{itemize}
\end{remark}

A main goal of this paper is to demonstrate how induction can be used to show that $T(\vec {\mathbf n},s,\vec {\mathbf 0})$ is true for many choices of $\vec {\mathbf n}$ and $s$. For this purpose it is enough to show that $$\dim \left[ T_{p_1}X^{\perp}\cap\ldots\cap
T_{p_s}X^{\perp}\right]$$ is less than or equal to the expected value for some choice of points $p_1,\dots, p_s$. By semicontinuity, establishing that the expected dimension holds in a particular case forces the expected dimension to hold in the general case. We reduce the size of a given problem through the specialization of sets of points. For instance, if $H\subseteq V_1$ is a subspace then we may specialize $t$ points among the points $p_1,\dots, p_s$
such that $p_i\in H$ for $i=1,\ldots, t$ then make our computation in this setting. If non-defectivity holds for a set of specialized points then it will hold for a set with the same number of general points.
This allows us to develop the following induction theorem.

\begin{thm0} (Subabundance Theorem)
\label{main1}

Let $n_1=n_1'+n_1''+1$, let $s=s'+s''$,
$a_2=a_2'+a_2''$, \dots, $a_k=a_k'+a_k''$. Suppose
\begin{itemize}

\item[(1)] $T(n_1', n_2, \dots, n_k;s' ; a_1+s'',a_2', \dots, a_k')$ is true and subabundant

\item[(2)] $T(n_1'', n_2, \dots n_k;s'' ; a_1+s',a_2'', \dots, a_k'')$ is true and subabundant

\end{itemize}
Then $T(n_1, \dots, n_k;s ; a_1,\dots, a_k)$ is true and subabundant.

\end{thm0}

\proof\hskip -3pt .   Let $H\subseteq V_1$ be a subspace of dimension $n_1'+1$ and let
$X'={\bf P}(H)\times{\bf P}(V_2)\times\dots \times{\bf P}(V_k)$
be embedded in ${\bf P}(H\otimes V_2\otimes\dots \otimes V_k)$. 
In the same way let $X''={\bf P}(V_1/H)\times{\bf P}(V_2)\times\dots \times {\bf P}(V_k)$
be embedded in ${\bf P}(V_1/H\otimes V_2\otimes\dots \otimes V_k)$. Consider $s$ points $p_1,\dots, p_s$ and specialize $p_i=v_{1,i}\otimes v_{2,i}\otimes \dots \otimes v_{k,i}$
in such a way that $v_{1,i}\in H$ for $i=1,\ldots ,s'$.  Let $f\colon V_1^{\vee}\otimes V_2^{\vee}\otimes \dots \otimes V_k^{\vee}\rig{}
H^{\vee}\otimes V_2^{\vee}\otimes\dots \otimes V_k^{\vee}$ be the natural projection.

By \lemref{restriction} we have
$f(T_{p_i}X^{\perp})=T_{p_i}X'^{\perp}$ for $i=1,\ldots ,s'$ and
$f(T_{p_i}X^{\perp})=G^1_{p_i}X^{\perp}$ for $i=s'+1,\ldots ,s$.
More precisely we have the exact sequences
$$0\rig{}G^1_{p_i}X''^{\perp}\rig{}T_{p_i}X^{\perp}\rig{}T_{p_i}X'^{\perp}\rig{}0$$
for $i=1,\ldots ,s'$ and the exact sequences
$$0\rig{}T_{[p_i]}X''^{\perp}\rig{}T_{p_i}X^{\perp}\rig{}G^1_{p_i}X^{\perp}\rig{}0$$
for $i=s'+1,\ldots ,s$ (where $[p_i]$ denotes the quotient class of $p_i$).

Combining these exact sequences yields
$$0\rig{}\cap_{i\le s'}G^1_{p_i}X''^{\perp}\bigcap\cap_{i>s'}
T_{[p_i]}X''^{\perp}\rig{}\cap_{i=1}^sT_{p_i}X^{\perp}\rig{}\cap_{i\le s'}T_{p_i}X'^{\perp}\bigcap\cap_{i>s'}G^1_{p_i}X^{\perp}.$$

We want to compute the dimension of the middle term $\cap_{i=1}^sT_{p_i}X^{\perp}$.
This explains why we have to include the spaces $G^i_{p_j}$ in the inductive procedure
from the very beginning.

Consider $a_1$ generic points $q_{1,1},\ldots ,q_{1,a_1}\in X$.
We get the exact sequences
$$0\rig{}G^1_{q_{1,i}}X''^{\perp}\rig{}G^1_{q_{1,i}}X^{\perp}\rig{}G^1_{q_{1,i}}X'^{\perp}\rig{}0$$
for $i=1,\ldots ,a_1$. 

Consider $a_2$ generic points $q_{2,1},\ldots ,q_{2,a_2}\in X$ and specialize
$q_{2,i}=v_{1,2,i}\otimes v_{2,2,i}\otimes \dots \otimes v_{k,2,i}$
in such a way that $v_{1,2,i}\in H$ for $i=1,\ldots ,a_2'$.

We get that
$$G^2_{q_{2,i}}X\simeq G^2_{q_{2,i}}X'$$
for $i=1,\ldots ,a_2'$ and that
$$G^2_{[q_{2,i}]}X''\simeq G^2_{q_{2,i}}X$$
for $i=a_2'+1,\ldots ,a_2$ (where $[q_{2,i}]$ denotes the quotient class of $q_{2,i}$).

In the same way, for  $a_t$ generic points $q_{t,1},\ldots ,q_{t,a_t}\in X$
we get that
$$G^t_{q_{t,i}}X\simeq G^t_{q_{t,i}}X'$$
for $i=1,\ldots ,a_t'$ and that
$$G^t_{[q_{t,i}]}X''\simeq G^t_{q_{t,i}}X$$
for $i=a_t'+1,\ldots ,a_t$.

Putting all of this together, we get the {\bf Fundamental Exact Sequence}
$$0\rig{}T_{s''}X''^{\perp}\cap G^1_{a_1+s'}X''^{\perp}\cap G^2_{a_2''}X''^{\perp}\cap\dots \cap G^k_{a_k''}X''^{\perp}
\rig{}$$
$$T_sX^{\perp}\cap G^1_{a_1}X^{\perp}\cap G^2_{a_2}X^{\perp}\cap\dots\cap G^k_{a_k}X^{\perp}\rig{}T_{s'}X'^{\perp}\cap G^1_{a_1+s''}X'^{\perp}\cap G^2_{a_2'}X'^{\perp}\cap\dots \cap G^k_{a_k'}X'^{\perp}.$$

By assumption (1),  the right term has codimension
$s'(1+n_1'+n_2+\dots +n_k)+(a_1+s'')(n_1'+1)+a_2'(n_2+1)+\dots +a_k'(n_k+1)$
in $H^{\vee}\otimes V_2^{\vee}\otimes \dots \otimes V_k^{\vee}$, meaning that all the intersections are transverse. By assumption (2), the left term 
has codimension
$s''(1+n_1''+n_2+\dots +n_k)+(a_1+s')(n_1''+1)+a_2''(n_2+1)+\dots +a_k''(n_k+1)$
in $(V_1/H)^{\vee}\otimes V_2^{\vee}\otimes \dots \otimes V_k^{\vee}$.
It follows that the middle term has codimension greater than or equal to
$s(1+\sum{n_i})+\sum (a_i)(n_i+1)$. Since this is the expected value, we have equality.\qedd

\

In the same way we have
\begin{thm0} (Superabundance Theorem)
\label{main2}

Let $n_1=n_1'+n_1''+1$, let $s=s'+s''$,
$a_2=a_2'+a_2''$, \dots, $a_k=a_k'+a_k''$. Suppose
\begin{itemize}

\item[(1)] $T(n_1', n_2, \dots, n_k;s' ; a_1+s'',a_2', \dots, a_k')$ is true and superabundant

\item[(2)] $T(n_1'', n_2, \dots n_k;s'' ; a_1+s',a_2'', \dots, a_k'')$ is true and superabundant

\end{itemize}
Then $T(n_1, \dots, n_k;s ; a_1,\dots, a_k)$ is true and superabundant.

\end{thm0}
\proof\hskip -3pt .   We proceed as in the previous theorem until we get to the {\bf Fundamental Exact Sequence}.
By assumption (1), the right term is zero. By assumption (2), the left term is zero.
It follows that the middle term is zero, as required.\qedd

\

\begin{coro0} If $T(n_1', n_2, \dots, n_k;s' ; a_1+s'',a_2',\dots,  a_k')$ and $T(n_1'', n_2, \dots, n_k;s'' ; a_1+s',a_2'',\dots, a_k'')$ are both true and equiabundant  then
$T(n_1, \dots, n_k;s ; a_1,\dots , a_k)$ is true and equiabundant.
\end{coro0}

\begin{remark} A simple but useful fact is that if $T(n_1,\dots,n_k;s;a_1,\dots,a_k)$ is true then $T(n_1,\dots,n_k,0;s;a_1,\dots,a_k,A)$ is true for any value of $A$.
\end{remark}

 It is important to note that if $n_1=1$ then we may take $n_1'=n_1''=0$ in \thmref{main1}
and \thmref{main2}. This allows us to reduce to a lower number of factors. Due to the importance of these cases, we state them explicitly as corollaries.

\begin{coro0} Let  $s=s'+s''$ and let
$a_j=a_j'+a_j''$, for $j=2,\ldots ,k$.  
\label{main5} 
Suppose that $(0,n_2,\dots,n_k;s;a_1,\dots,a_k)$ is subabundant then $T(0,n_2,\dots,n_k;s;0,a_2,\dots,a_k)$ is true if and only if  $T(0,n_2,\dots,n_k;s;a_1,a_2,\dots,a_k)$ is true.
\end{coro0}

\proof\hskip -3pt .   We reduce to Theorem \ref{main1} because the corresponding condition $G^{1{\perp}}$ is of codimension one and is independent from the other conditions provided subabundancy is satisfied. If $a_1$ is such that $(1, n_2, \ldots,  n_k;s ; a_1, \ldots,  a_k)$ is superabundant then $T(1, n_2, \ldots,  n_k;s ; a_1, \ldots,  a_k)$ is also true as the ambient space is filled.\qedd

\

\begin{coro0}
\label{main6}
Let  $s=s'+s''$ and let
$a_j=a_j'+a_j''$, for $j=2,\ldots ,k$.  If both $T( n_2,\ldots,  n_k;s' ; a_2',\ldots,  a_k')$ and $T( n_2, \ldots,  n_k;s'' ; a_2'',\ldots,  a_k'')$ are true and superabundant then $T(1, n_2, \ldots,  n_k;s ; a_1,a_2, \ldots,  a_k)$ is true (and superabundant).
\end{coro0}

\begin{remark} Theorem \ref{main1} and Theorem \ref{main2} should be viewed as a generalization of the Splitting Method of B\"urgisser, Clausen and Shokrollahi  from the case of 3 factors to the case of $k$ factors \cite{BCS}. The proof given in the present paper takes a more geometric and homological point of view and mirrors the ideas of Alexander-Hirschowitz and Terracini in the use of degeneration arguments\cite{AH,T}. A proof written purely in the language of tensor algebra would also be natural following the approach of B\"urgisser et al. This would have the advantage of conciseness but the geometry would be pushed more to the background.

Recall that if $X\subseteq \P m$ is a variety and if $W_1,\dots, W_p$ are subvarieties of $X$ then $J(W_1,\dots, W_p)$ is called a partial secant variety to $X$. In the particular case when $X$ is the Segre variety $\P {{n_1}}\times \dots \times \P {{n_k}}$,
the linear space $L$ spanned by $T_sX+G^1_{a_1}X+G^2_{a_2}X+\dots +G^k_{a_k}X\subseteq T_{(1,\ldots ,1)}$ should be seen as the tangent space to a particular partial secant variety of $X$. The expression $T_sX$ corresponds to computing the tangent space to $X$ at $s$ general points. The expression $G^1_pX$ corresponds to computing the tangent space at a general point, $p=v_1\times \dots \times v_k$, of a subvariety of $X$ of the form $\P {{n_1}} \times v_2 \times \dots \times v_k$. Such a subvariety is a $\P {{n_1}}$ sitting inside $X$. The expression $G^1_{a_1}X$ corresponds to computing the span of the tangent spaces to $a_1$ different such subvarieties for $a_1$ different choices of $p$.  Similarly, each of the other $G^i_{a_i}$ represent the span of tangent spaces to $a_i$ different varieties in the family of $\P {{n_i}}$'s obtained by fixing all but the $i^{th}$ coordinate. Thus  viewing $L$ as a tangent space at a general point of the join of a collection of $s+a_1+\dots +a_k$ subvarieties of $X$ follows as an immediate application of Terracini's Lemma as stated in \cite{A}. Furthermore, $a_1+\dots +a_k$ of the subvarieties are linear spaces inside $X$.
\end{remark}

Theorem \ref{main1} should be viewed as a way of computing of the dimension of a secant variety by applying semicontinuity arguments to the computation of the dimension of smaller partial secant varieties arising from specializations of points.
It is clear that after a finite number of applications of the previous two theorems, we may reduce ourselves to the four 
projective varieties
$$\P 1\times\P 1\times\P 1,\qquad  \P 1\times\P 1\times\P 2,\qquad \P 1\times\P 2\times\P 2
\hskip 10pt {\rm and}\hskip 10pt  \P 2\times\P 2\times\P 2.$$
The importance of this reduction is emphasized in the following proposition,
which was essentially proved by Strassen:
\begin{prop0} \label{strass} \cite{S} Suppose $T(\vec {\mathbf n},s,\vec {\mathbf a})$ is true.
\begin{itemize}
\item[(i)] If $T(\vec {\mathbf n},s,\vec {\mathbf a})$ is  subabundant
 and if $\vec {\mathbf n}'\geq \vec {\mathbf n}$ then $T(\vec {\mathbf n}',s,\vec {\mathbf a})$ is true and subabundant.
\item[(ii)] If $T(\vec {\mathbf n},s,\vec {\mathbf a})$ is superabundant
 and if $\vec {\mathbf n}'\leq \vec {\mathbf n}$ then $T(\vec {\mathbf n}',s,\vec {\mathbf a})$ is true and superabundant.
\end{itemize}
\end{prop0}

\proof\hskip -3pt .   In order to prove the first statement we can reduce to 
the case where $n_i=n'_i$ for $i=1,\ldots ,k-1$ and  $n_k+1=n'_k$.
Fix a splitting $V'_k=V_k\oplus <v>$. This induces an inclusion $X=\P {{n_1}}\times\ldots\times \P {{n_k}}
\subset \P {{n'_1}}\times\ldots\times \P {{n'_k}}=X'$ corresponding to the splitting
\begin{equation}\label{vsplit}
V'_1\otimes\ldots\otimes V'_k=(V_1\otimes\ldots\otimes V_k)\hskip 6pt\oplus\hskip 6pt (V_1\otimes\ldots\otimes V_{k-1}\otimes<v>)
\end{equation}
 Pick a point $p=v_1\otimes\ldots\otimes v_k
\in \P {{n_1}}\times\ldots\times \P {{n_k}}$.
In affine notation we have $$T_pX'=T_pX\oplus <v_1\otimes\ldots\otimes v_{k-1}\otimes v>$$
$$G^i_pX'=G^i_pX\hbox{\ for\ }i=1,\ldots ,k-1$$  
$$G^k_pX'=G^k_pX\oplus <v_1\otimes\ldots\otimes v_{k-1}\otimes v>$$
and these splitting are compatible with (\ref{vsplit}).  Now it is easy to check
that if $T_{p_i}X$ and $G^t_{q_{t,i}}X$ are transversal then 
$T_{p_i}X'$ and $G^t_{q_{t,i}}X'$ are also transversal.
The second statement proceeds in an analogous manner.\qedd

\

\begin{remark} 
We can utilize Proposition~\ref{strass} for higher dimensional Segre varieties by ``{\rm padding with zeroes}".
 For instance, if $T(n_1,n_2,n_3;s;0,0,0)$ is true and subabundant then we can pad with a zero to obtain that
 $T(n_1,n_2,n_3,0;s;0,0,0,0)$ is true and subabundant. As a consequence, by Proposition~\ref{strass},
 $T(n_1,n_2,n_3,n_4;s;0,0,0,0)$ is true and subabundant for any $n_4$ (being sure to keep $s$ fixed).
\end{remark}

\begin{notation} We introduce the notation $b*T(\vec{{\bf n}};s;\vec{\bf {a}})$ to denote $b$ identical statements of the form $T(\vec{{\bf n}};s;\vec{\bf {a}})$.
\end{notation}

For the four projective varieties $\hskip 2pt\P 1\times\P 1\times\P 1, \hskip 5pt \P 1\times\P 1\times\P 2, \hskip 5pt \P 1\times\P 2\times\P 2
\hskip 5pt {\rm and}\hskip 5pt  \P 2\times\P 2\times\P 2,$
we list the $4$-tuples $(s;a_1,a_2,a_3)$ where the statement $T(n_1,n_2,n_3;s;a_1,a_2,a_3)$ is not true.
For the varieties $\P 1\times\P 2\times\P 2
\hskip 5pt {\rm and}\hskip 5pt  \P 2\times\P 2\times\P 2$, we divide the list into the minimal cases and the non-minimal cases. The defectivity of each of the non-minimal cases follows directly from the defectivity of one of the minimal cases. The defectivity of the minimal cases are all established by the elementary arguments given in the following 3 lemmas. The non-defectivity of the cases not appearing on these lists can be established by explicit computation.

\begin{lemma0}
Let $(s,a_2)=(1,0)$ or $(0,1)$. 
Suppose that the following inequality holds: 
\[
a_3+sn_1+n_2 \geq (n_1+1)(n_2+1).  
\]
Then $T(n_1,n_2,n_3;s;0,a_2,a_3)$ is false.  
\end{lemma0}
\proof\hskip -3pt . 
Let $X=\mathbf{P}^{n_1} \times \mathbf{P}^{n_2} \times \mathbf{P}^{n_3}$ 
and let $q_1, \dots, q_{a_3}$ be general points of $X$. 
Note that $X$ can be viewed as a 
$(\mathbf{P}^{n_1} \times \mathbf{P}^{n_2})$-fibration over $\mathbf{P}^{n_3}$. 

Suppose that $(s,a_2)=(1,0)$. 
Given a point $p$ of $X$, there is a fiber 
$\mathbf{P}^{n_1} \times \mathbf{P}^{n_2} \subset \mathbf{P}^{(n_1+1)(n_2+1)-1}$, 
denoted $Q$, which contains $p$.  
For each $i \in \{1, \dots, a_3\}$, the projectivization of $G_{q_i}^3 X$ is   
a horizontal $n_3$-plane, which meets $Q$ in a single point. 
The $a_3$ points as obtained above span a $\mathbf{P}^{a_3-1} 
\subset \mathbf{P}^{(n_1+1)(n_2+1)-1}$, 
and  if  $\dim \mathbf{P}^{a_3-1}+ \dim Q \geq \mathbf{P}^{(n_1+1)(n_2+1)-1}$, 
then every tangent space to $Q$, and thus the 
tangent space to $Q$ at $p$,  must intersect $\mathbf{P}^{a_3-1}$. 
Therefore, $\mathbf{P}\left(G_{a_3}^3X \right)$ meets the tangent space to $X$ 
at $p$. This implies that $T(n_1,n_2,n_3;1;0,0,a_3)$ fails. 

In a similar way, one can show that $T(n_1,n_2,n_3;0;0,1,a_3)$ fails. 
Given a point $p$, 
there is a fiber 
$\mathbf{P}^{n_1} \times \mathbf{P}^{n_2} \subset \mathbf{P}^{(n_1+1)(n_2+1)-1}$, 
which contains the $n_2$-plane $\mathbf{P}(G_p^2X)$.   
This $n_2$-plane $\mathbf{P}(G_p^2X)$ and 
the $(a_3-1)$-plane $\mathbf{P}^{a_3-1}$ as obtained above must intersect, 
because we have 
$\dim \mathbf{P}^{a_3-1}+ \dim \mathbf{P}^{n_2} \geq \mathbf{P}^{(n_1+1)(n_2+1)-1}$ 
by the inequality as given above. 
In particular, $\mathbf{P}\left(G_{a_3}^3X \right)$ meets  $\mathbf{P}(G_p^2X)$, 
which implies that $T(n_1,n_2,n_3;0;0,1,a_3)$ fails. 
\qedd

\

\begin{lemma0}
$T(\vec{{\bf n}};s;\vec{\bf {a}})$ is false for the cases
$(\vec{{\bf n}};s;\vec{{\bf a}})=(1,2^2;2;0,0,2)$, $(2^2;2;0,0,4)$ and  
$(2^3;3;0,1,1)$. 
\end{lemma0} 
\proof\hskip -3pt . 
The main idea of this lemma is to use the contrapositive of \thmref{main1}. 
Note that $(1,2,5;4;0,0,0)$ is unbalanced (see Lemma \ref{unbalanced} and Definition \ref{defunbalanced}). 
Thus the statement  $T(1,2,5;4;0,0,0)$ is equiabundant, but not true.  
One can reduce this statement to the equiabundant statement $2*T(1,2,2;2,0,2)$. 
So the fact that $T(1,2,5;4;0,0,0)$ is not true implies that $T(1,2,2;2,0,2)$ is not true.  

In a similar manner, we can prove that $T(2,2,2;2,0,4)$ is not true. 
Note that $(2,2,8;6;0,0,0)$ is unbalanced (see Proposition 4.1) and $T(2,2,8;6;0,0,0)$ is false. 
One can reduce this statement to $3*T(2,2,2;2;0,0,4)$. 
So the fact that $T(2,2,8;6;0,0,0)$ is not true implies 
that $T(2,2,2;2;0,0,4)$ is not true.  

By \propref{p2p3p3}, the subabundant statement $T(2,3,3;5;0,0,0)$ is false. 
This implies that one of the statements 
$T(2,2,3;4;0,1,0)$ and $T(2,0,3;1;0,4,0)$ is false. 
Clearly the second statement is true, and so $T(2,2,3;4;0,1,0)$ cannot be true. 
Since the $T(2,2,3;4;0,1,0)$ can be reduced to the subabundant 
statements $T(2,2,2;3;0,1,1)$ or $T(2,2,0;1,0,0,3)$, 
we can say that either $T(2,2,2;3;0,1,1)$ or $T(2,2,0;1,0,0,3)$ is false. 
Since the second statement is true, we can conclude that 
$T(2,2,2;3;0,1,1)$ is false, which completes the proof. 
\qedd

\

\begin{lemma0} $T(2,2,2;4;0,0,0)$ is false.
\end{lemma0}
\proof\hskip -3pt .  This case is well known. The geometrical explanation 
is the following.
Given four points in $X=\P 2\times\P 2\times\P 2\subset \P {{26}}$,
we can project on each factor, and get isomorphisms that identify the three factors. The diagonal surface, after this identification, is the $3$-Veronese embedding of $\P 2$,
which contains the four original points and span a linear $\P 9$. The four tangent spaces
to $X$ at these points meet the $\P 9$ in dimension $\ge 2$, and
 the dimension of $\sigma_4(X)$ is at most $9+4\cdot 4=25$. 
\qedd

\

\begin{prop0} \label{listdef} The following is a complete list of the defective $(\vec {\mathbf n},s,\vec {\mathbf a})$  with $\vec {\mathbf n}=(n_1,n_2,n_3)$ and $1\leq n_1,n_2,n_3\leq 2$. The list is given as $(s;a_1,a_2,a_3)$.
\begin{itemize}
        \item[(i)] {$\P 1\times\P 1\times\P 1$}\qquad
Up to permutation of the three factors the list is

Minimal: $(0;0,1,3)$, $(1;0,0,2)$ 

        \item[(ii)] {$\P 1\times\P 1\times\P 2$}\qquad
        Up to permutation of the first two factors the list is
        
Minimal: $(0;0,1,3)$, $(0;0,4,1)$, $(0;1,5,0)$, $(1;0,3,0)$,     $(1;0,0,2)$. 

        \item[(iii)] {$\P 1\times\P 2\times\P 2$}\qquad
Up to permutation of the last two factors the list is
        
Minimal: $(0;0,1,4),  (0;7,0,1), (0;1,0,5), (1;0,0,3), (1;5,0,0), (2;0,0,2)$

Non-minimal: $(1;6,0,0), (0;0,1,5), (0;0,2,4), (0;1,1,4), (1;0,0,3),$ 

$\hskip 70pt (1;0,0,4), (1;0,1,3), (1;1,0,3)$.

\item[(iv)] {$\P 2\times\P 2\times\P 2$}\qquad
Up to permutation of the three factors the list is
        
Minimal: $(0;0,1,7),   (1;0,0,5), (2;0,0,4), (3;0,1,1), (4;0,0,0)$

Non-minimal: $(0;1,1,7), (0;0,2,7), (0;0,1,8), (1;0,0,6), (1;0,1,5)$

\end{itemize}
\end{prop0}

\proof\hskip -3pt . The defectivity of the minimal cases follow from the previous 3 lemmas. The non-minimal cases follow from the minimal cases. The non-defectivity of the $(\vec {\mathbf n},s,\vec {\mathbf a})$ not appearing on the list can be shown by explicit computation. \qedd

\

We will now illustrate the inductive method of Theorem~\ref{main1} and Theorem~\ref{main2} in a series of examples. The strategy is to reduce a problem involving a more complicated variety to known cases on simpler varieties. By Remark~\ref{usefulrem}, in order to establish the non-defectivity of all secant varieties to a given Segre variety, it is enough to check the truth of statement $T(\vec {\mathbf n},s,\vec {\mathbf 0})$ for the largest $s$ for which $(\vec {\mathbf n},s,\vec {\mathbf 0})$ is subabundant and for the smallest  $s$ for which $(\vec {\mathbf n},s,\vec {\mathbf 0})$ is superabundant.

\begin{exa} In this example we show that
$X=\P 3\times\P 3\times\P 3$ has no defective secant varieties
(already known by Lickteig). This is reduced to showing that $\dim\sigma_6(X)=59$ and that $\sigma_7(X)$ fills the ambient space.

\

In order to prove that $\dim\sigma_6(X)=59$, we need to establish $T(3,3,3;6;0,0,0)$.
We have $$T(1,3,3;3;3,0,0) \hskip 10pt {\rm and} \hskip 10pt T(1,3,3;3;3,0,0) \Rightarrow T(3,3,3;6;0,0,0)$$
$$T(1,1,3;2;1,1,0)\hskip 10pt {\rm and} \hskip 10pt  T(1,1,3;1;2,2,0) \Rightarrow T(1,3,3;3;3,0,0))$$
$$T(1,1,1;1;1,0,1)\hskip 10pt {\rm and} \hskip 10pt  T(1,1,1;1;0,1,1) \Rightarrow T(1,1,3;2;1,1,0)$$
$$T(1,1,1;1;1,1,0)\hskip 10pt {\rm and} \hskip 10pt  T(1,1,1;0;1,1,1) \Rightarrow T(1,1,3;1;2,2,0)$$
But, $T(1,1,1;1;1,0,1), T(1,1,1;1;0,1,1),T(1,1,1;1;1,1,0)$ and $T(1,1,1;0;1,1,1)$ are all true, thus $T(3,3,3;6;0,0,0)$ is true and $\dim\sigma_6(X)=59$.

\

In order to prove that $\sigma_7(X)$ fills the ambient space we 
need $T(3,3,3;7;0,0,0)$ to be true.
We have $$T(1,3,3;4;3,0,0) \hskip 10pt {\rm and} \hskip 10pt T(1,3,3;3;4,0,0) \Rightarrow T(3,3,3;7;0,0,0)$$
$$T(1,1,3;2;1,2,0)\hskip 10pt {\rm and} \hskip 10pt  T(1,1,3;2;2,2,0) \Rightarrow T(1,3,3;4;3,0,))$$
$$T(1,1,3;2;1,1,0)\hskip 10pt {\rm and} \hskip 10pt  T(1,1,3;1;3,2,0) \Rightarrow T(1,3,3;3;4,0,0)$$
$$T(1,1,1;1;1,1,1)\hskip 10pt {\rm and} \hskip 10pt  T(1,1,1;1;0,1,1) \Rightarrow T(1,1,3;2;1,2,0)$$
$$T(1,1,1;1;1,1,1)\hskip 10pt {\rm and} \hskip 10pt  T(1,1,1;1;1,1,1) \Rightarrow T(1,1,3;2;2,2,0)$$
$$T(1,1,1;1;1,0,1)\hskip 10pt {\rm and} \hskip 10pt  T(1,1,1;1;0,1,1) \Rightarrow T(1,1,3;2;1,1,0)$$
$$T(1,1,1;1;1,1,0)\hskip 10pt {\rm and} \hskip 10pt  T(1,1,1;0;2,1,1) \Rightarrow T(1,1,3;1;3,2,0)$$
The proof follows from the last 4 implications, thus $T(3,3,3;7;0,0,0)$ is true and $\sigma_7(X)$ fills the ambient space.
\end{exa}

\begin{exa}  In this example we show that
$X=\P 5\times\P 5\times\P 5$ has no defective secant varieties
(already known by Lickteig). This is reduced to showing that $\dim\sigma_{13}(X)=207$ and that $\sigma_{14}(X)$ fills the ambient space.
The example is shown in some detail to emphasize that the strategy of reduction can be tricky.

\

In order to prove that $\dim\sigma_{13}(X)=207$, we need to establish that $T(5,5,5;13;0,0,0)$ is true.
If we use \thmref{main1} to reduce to $T(2,5,5;7;6,0,0), T(2,5,5;6;7,0,0)$ then we find that the 7-tuple $(2,5,5;7;6,0,0)$ is not subabundant!

We modify our strategy
and  reduce to
$$T(1,5,5;4;9,0,0)\qquad T(3,5,5;9;4,0,0).$$
Then $T(1,5,5;4;9,0,0)$ can reduce to 
$$T(1,2,5;2;5,2,0)\quad T(1,2,5;2;4,2,0).$$ Since $T(1,2,5;2;5,2,0) \Rightarrow T(1,2,5;2;4,2,0)$ (see Remark \ref{usefulrem} (ii)), it is enough to consider
$T(1,2,5;2;5,2,0)$. This reduces to
$$T(1,2,2;1;3,1,1)\quad T(1,2,2;1;2,1,1).$$
Both of these statements are true.

Now we reduce $T(3,5,5;9;4,0,0)$ to 
$$T(3,3,5;6;3,3,0)\quad T(3,1,5;3;1,6,0)$$
which reduce respectively to

$$(1)\qquad T(3,3,2;3;1,2,3)\quad T(3,3,2;3;2,1,3)$$
$$(2)\qquad T(3,1,2;2;0,3,1)\quad T(3,1,2;1;1,3,2).$$
(1) consists of two equivalent cases. We reduce $T(3,3,2;3;1,2,3)$ to 
$$T(1,3,2;1;3,2,1)\quad T(1,3,2;2;2,0,2)$$
and finally to
$$T(1,1,2;0;3,3,0)\quad T(1,1,2;1;0,2,1)$$
$$T(1,1,2;1;1,1,1)\quad T(1,1,2;1;1,1,1).$$
These last four statements are true.

(2) reduces to 
$$T(1,1,2;1;1,1,1)\quad T(1,1,2;1;1,2,0)$$ $$T(1,1,2;1;1,1,1)\quad T(1,1,2;0;2,2,1).$$
These last four statements are true.

Thus we have proved that $\dim\sigma_{13}(X)=207$.

\

In order to prove that $\sigma_{14}(X)$ fills the ambient space,
we  reduce by \thmref{main2} to
$$T(2,5,5;7;7,0,0)\qquad T(2,5,5;7;7,0,0).$$
Then $T(2,5,5;7;7,0,0)$ reduces to
$$T(2,2,5;4;2,3,0)\qquad T(2,2,5;3;5,4,0)$$
which reduces to
$$T(2,2,2;3;0,1,1)\qquad T(2,2,2;1;2,2,3)$$
$$T(2,2,2;2;2,2,1)\qquad T(2,2,2;1;3,2,2).$$

Unfortunately the statement $T(2,2,2;3;0,1,1)$ is not true,
so we have not proven anything.

We change our strategy and from 
$T(2,5,5;7;7,0,0)$
we reduce to
$$T(2,1,5;3;1,4,0)\qquad T(2,3,5;4;6,3,0).$$
Then $T(2,1,5;3;1,4,0)$ reduces to
$$T(2,1,2;2;0,2,1)\qquad T(2,1,2;1;1,2,2)$$
while $T(2,3,5;4;6,3,0)$ reduces to 
$$T(2,1,5;2;3,5,0)\qquad T(2,1,5;2;3,5,0)$$
and finally to
$$T(2,1,2;1;2,2,1)\qquad T(2,1,2;1;1,3,1).$$

Now all the final reduced statements are true and we have proved that
$\sigma_{14}(X)$  fills the ambient space. 
\end{exa}

Let us now show an example which seems to be new.

\begin{exa} Consider $X=\P 4\times\P 4\times\P 7\subset\P {{199}}$. We have
$\lfloor 200/16\rfloor=12$, $\lceil 200/16\rceil=13$.
In order to show that $\sigma_{12}(X)$ has the expected dimension $191$,
we reduce $T(4,4,7;12;0,0,0)$ by \thmref{main1} to
$$T(2,4,7;7;5,0,0)\qquad T(1,4,7;5;7,0,0)$$
The first one reduces to
$$T(2,2,7;4;3,3,0)\qquad T(2,1,7;3;2,4,0)$$
and the second one reduces to
$$T(1,2,7;3;4,2,0)\qquad T(1,1,7;2;3,3,0)$$
These last four statements reduce respectively to
$$(1)\qquad T(2,2,1;1;1,0,3)\quad T(2,2,1;1;0,1,3)\quad T(2,2,1;1;1,1,3)\quad T(2,2,1;1;1,1,3)$$
$$(2)\qquad T(2,1,1;1;0,1,2)\quad T(2,1,1;1;1,0,2)\quad T(2,1,1;1;1,0,2)\quad T(2,1,1;0;0,3,3)$$
$$(3)\qquad T(1,2,1;1;1,0,2)\quad T(1,2,1;1;0,1,2)\quad T(1,2,1;1;0,1,2)\quad T(1,2,1;0;3,0,3)$$
$$(4)\qquad T(1,1,1;1;1,0,1)\quad T(1,1,1;1;0,1,1)\quad T(1,1,1;0;1,1,2)\quad T(1,1,1;0;1,1,2).$$
These statements are all true and we conclude that $\dim\sigma_{12}(X)=191$.

\

To show that $\sigma_{13}(X)$ fills $\P{{199}}$,
we reduce $T(4,4,7;13;0,0,0)$ by \thmref{main2} to
$$T(2,4,7;8;5,0,0)\qquad T(1,4,7;5;8,0,0).$$
The first one reduces to
$$T(2,2,7;5;3,3,0)\qquad T(2,1,7;3;2,5,0)$$
and the second one reduces to
$$T(1,2,7;3;5,2,0)\qquad T(1,1,7;2;3,3,0).$$
These last four statements reduce respectively to
$$(1)\qquad T(2,2,1;2;0,0,3)\quad T(2,2,1;1;1,1,4)\quad T(2,2,1;1;1,1,4)\quad T(2,2,1;1;1,1,4)$$
$$(2)\qquad T(2,1,1;1;0,2,2)\quad T(2,1,1;1;1,0,2)\quad T(2,1,1;1;1,0,2)\quad T(2,1,1;0;0,3,3)$$
$$(3)\qquad T(1,2,1;1;2,0,2)\quad T(1,2,1;1;0,1,2)\quad T(1,2,1;1;0,1,2)\quad T(1,2,1;0;3,0,3)$$
$$(4)\qquad T(1,1,1;1;1,0,1)\quad T(1,1,1;1;0,1,1)\quad T(1,1,1;0;1,1,2)\quad T(1,1,1;0;1,1,2).$$
These last statements are all true and we conclude that $\sigma_{13}(X)$ fills the ambient space.

\end{exa}

\section{Classification of Segre varieties with defective $r$-secant varieties, $r\le 6$}

In this section, $X=\P {{n_1}}\times\ldots\times\P {{n_k}}$ with
$k\ge 3$ and $n_1\le\ldots \le n_k$. We classify Segre varieties, $X$, for which $\sigma_r(X)$ is defective with $r\leq 6$.  We recall that no Segre variety with 3 or more factors has a defective 2-secant variety. 

Following \cite{BCS}, the typical tensor rank of a format $(n_1,\ldots ,n_k)$ is the smallest integer
$s$ such that $\sigma_s(\P {{n_1}}\times\ldots\times \P {{n_k}} )$ fills the ambient space, and it is denoted by
$\underline R (n_1,\ldots ,n_k)$. Equivalently, the generic tensor in $V_1\otimes\ldots\otimes V_k$
where $\dim V_i=n_i+1$ is the sum of  $\underline R (n_1,\ldots ,n_k)$ (and not less) tensors of rank one.
We use the projective notation, so that our  $\underline R (n_1,\ldots ,n_k)$ corresponds to
$\underline R (n_1+1,\ldots ,n_k+1)$  of \cite{BCS}.
Obviously we have
$$\left\lceil \frac{\prod(n_i+1)}{1+\sum n_i} \right\rceil\le \underline R (n_1,\ldots ,n_k)$$
and in particular
$$\left\lceil \frac{(n+1)^k}{nk+1} \right\rceil\le \underline R (n^k).$$

\

The following lemma is well-known (see [CGG1, Proposition 3.3]). 
\begin{lemma0} \label{unbalanced}
Let $X=\mathbf{P}^{n_1} \times \dots \times
\mathbf{P}^{n_k}$, $1 \leq n_1 \leq  \dots \leq n_k$. 
Suppose that 
\[
\prod_{i=1}^{k-1} (n_i+1) - \sum_{i=1}^{k-1}  n_i < d 
< \min\left\{ \prod_{i=1}^{k-1} (n_i+1), n_k+1\right\}. 
\]
Then $X$ has a defective $d$-secant variety. 
\end{lemma0}
\proof\hskip -3pt .  Pick $d$ general points on $X$ where $d$ satisfies the conditions of the Lemma. 
Since $d < n_k+1$,  there exists a subvariety 
$V=\mathbf{P}^{n_1} \times \dots \times \mathbf{P}^{n_{k-1}} \times \mathbf{P}^{d-1} \subseteq X$, which contains these $d$ points. 
Let $N(d)=d\prod_{i=1}^{k-1} (n_i+1) -1$ and $N=\prod_{i=1}^k (n_i+1) -1$. 
The span of $V$ is $\mathbf{P}^{N(d)}\subseteq \mathbf{P}^N$. 
Thus, the linear subspace spanned by the tangent spaces of $X$ at the $d$ points 
has dimension at most $F(d)-1$, 
where $F(d)=d\left[\prod_{i=1}^{k-1}(n_i+1)+(n_k+1-d)\right]$. 
Then, by the assumption as given above, we have 
\begin{eqnarray*} 
d\left(\sum_{i=1}^k n_i+1\right)-F(d) 
&=& 
d\left(\sum_{i=1}^k n_i+1\right)-d\left[\prod_{i=1}^{k-1}(n_i+1)+(n_k+1-d)\right] \\
&=&
d\left[\sum_{i=1}^{k-1} n_i - \prod_{i=1}^{k-1}(n_i+1) +d \right] >0
\end{eqnarray*}
and 
\begin{eqnarray*}
\prod_{i=1}^k (n_i+1) -F(d)
&=& 
d^2-d\left[ \prod_{i=1}^{k-1}(n_i+1)+(n_k+1)\right]+\prod_{i=1}^k (n_i+1) \\
&=& 
\left[d-\prod_{i=1}^{k-1}(n_i+1)\right]\left[d-(n_k+1)\right]>0. 
\end{eqnarray*}
So $F(d) < 
\min\left\{ d\left(\sum_{i=1}^k n_i+1\right), \prod_{i=1}^k (n_i+1) \right\}$.
An application of Terracini's lemma shows that $X$ has a defective $d$-secant variety.    
\qedd

\

\begin{defn0} \label{defunbalanced}  Suppose $\vec {\mathbf n}=(n_1,\dots, n_k)$ with $n_1\le \dots \le n_k$.
\begin{itemize}
\item  $\vec {\mathbf n}$ is called {\bf balanced} if $n_k\le\prod_{i=1}^{k-1}(n_i+1)- \sum_{i=1}^{k-1}n_i$. 
\item  $\vec {\mathbf n}$ is called {\bf unbalanced} if $n_k-1\ge\prod_{i=1}^{k-1}(n_i+1)- \sum_{i=1}^{k-1}n_i$.
\end{itemize}
\end{defn0}

Thus Lemma \ref{unbalanced} states that if $\vec {\mathbf n}=(n_1,\dots, n_k)$ is unbalanced then ${\mathbf P}^{\vec {\mathbf n}}$ is defective.
The following proposition is often useful.

\begin{prop0}\label{balancuseful}
Let $\vec {\mathbf n}=(n_1,\dots, n_k)$ be balanced.
If $s\le n_k$ then $T( \vec {\mathbf n},s,0^k)$ is true and subabundant.
\end{prop0}

\proof\hskip -3pt .  It is sufficient to check the statement for $s=n_k$.
By assumption we have 
$$\sum_{i=1}^{k}  n_i\le  \prod_{i=1}^{k-1} (n_i+1)$$
After multiplying by $(n_k+1)$ we obtain
$$(1+\sum_{i=1}^{k}  n_i)n_k\le \sum_{i=1}^{k-1}  n_i+(1+\sum_{i=1}^{k}  n_i)n_k\le  \prod_{i=1}^{k} (n_i+1).$$
This implies that $( \vec {\mathbf n},n_k,0^k)$ is subabundant. By \thmref{main1}, $( \vec {\mathbf n},n_k,0^k)$
reduces to $T(n_1,\dots, n_{k-1},0;0,0^{k-1},n_k)$
and  $n_k*T(n_1,\dots, n_{k-1},0;1,0^{k-1},n_k-1)$.
Since both of these statements are true, we are done.\qedd

\

The following theorem sets completely the defective behaviour of higher secant varieties in the unbalanced cases,
and completes the Prop 3.3. in \cite{CGG1}. This has also been observed as part of Theorem 2.4 in \cite{CGG4}.

\begin{thm0}
Let $\vec {\mathbf n}=(n_1,\dots, n_k)$ be unbalanced.
\begin{itemize}
\item[(i)] $T( \vec {\mathbf n},s,0^k)$ is true and subabundant if and only if 
$s\le \prod_{i=1}^{k-1} (n_i+1) - \sum_{i=1}^{k-1}  n_i$.
\item[(ii)] $\underline R(\vec {\mathbf n}) =\min\{n_k+1,\hskip 2pt\prod_{i=1}^{k-1} (n_i+1)\}$
\end{itemize}
\end{thm0}

\proof\hskip -3pt .  The  ``only if" part of (i) is \lemref{unbalanced}. In order to prove the ``if" part,
set $n_k'= \prod_{i=1}^{k-1} (n_i+1) - \sum_{i=1}^{k-1}  n_i $.
It is enough to check that $T(\vec {\mathbf n};n_k';0^k)$ is true and subabundant.
By assumption we have $n_k'\le n_k-1$, moreover
$(n_1,\dots, n_{k-1},n_k')$ is balanced. By  \propref{balancuseful}
$T( n_1,\dots, n_{k-1},n_k';n_k';0^k)$  is true and subabundant.
The thesis follows by \propref{strass}.

Statement (ii) follows from Theorem 3.1 in \cite{CGG1}. \qedd

\

\begin{thm0}
$\sigma_3(\P {{n_1}}\times\ldots\times\P {{n_k}})$  is non-defective with the following exceptions:
 
$(n_1, n_2, n_3)=(1,1,a)$ with $a\ge 3$

$(n_1, n_2, n_3,n_4)=(1,1,1,1).$ 
\end{thm0}

\proof\hskip -3pt .   First we prove the theorem for $k=3$:  Since $T(1,2,2;3;0,0,0)$ is true and subabundant, from Proposition~\ref{strass}, we know that  $\sigma_3(X)$ has the expected dimension if $n_1\geq1, n_2\geq 2, n_3\geq2$. Hence, we may assume
$n_1=n_2=1$. $T(1,1,a;3;0,0,0)$ is true for $a=1,2$. $T(1,1,a;3;0,0,0)$ is false for $a\ge 3$ by Lemma \ref{unbalanced}. 

To prove the theorem for $k\ge 4$, it is enough to exhibit three points such that their tangent spaces are independent.
It is known that $\dim\sigma_3( \P 1\times\P 1\times \P 1\times\P 1)$ is smaller than expected so with four factors assume that $n_4\ge 2$. Then choose $(e_0,e_0,e_0,e_0)$, $(e_1,e_1,e_1,e_1)$, $(e_0+e_1,e_0,e_1,e_2)$.
With at least five factors choose $(e_0,e_0,e_0,e_0,e_0,*)$, $(e_1,e_1,e_1,e_0,e_0,*)$, $(e_0,e_0,e_1,e_1,e_1,*)$.
\qedd

\

\begin{thm0}
$\sigma_4(\P {{n_1}}\times\ldots\times\P {{n_k}})$  is non-defective with the following exceptions:
 
$(n_1, n_2, n_3)=(1,2,a)$ with $a\ge 4$

$(n_1, n_2, n_3)=(2,2,2)$.
\end{thm0}

\proof\hskip -3pt .   It is known that $T(1,1,1,1,1;4;0,0,0,0,0)$ is true. Thus there are no exceptions with $k\ge 5$. To treat the case $k=4$ we consider the equiabundant case $(1,1,1,2;4;0,0,0,0)$ is true. By \thmref{main1}, $T(1,1,1,2;4;0,0,0,0)$ reduces to twice $T(0,1,1,2;2;2,0,0,0)$. Since this is known to be true, there are no exceptions with $k=4$.

To treat the case $k=3$, we start with the known fact that $\dim\sigma_4(\P 2\times\P 2\times\P 2)$ is smaller than expected. So let us begin by proving
that $T(2,2,3;4;0,0,0)$ is true (and subabundant).
Indeed we reduce by \thmref{main1} to
$T(2,2,1;2;0,0,2)$ which is true.
Hence if $n_1\ge 2$ the theorem holds and we may assume $n_1=1$.

Let us now prove that  $T(1,3,3;4;0,0,0)$ is true (and subabundant).
We reduce by \thmref{main1} to
twice $T(1,1,3;2;0,2,0)$ and then reduce to four $T(1,1,1;1;0,1,1)$.
This is known to be true, hence if $n_2\ge 3$ the theorem holds and we may assume $n_2=2$. 

 $T(1,2,a;4;0,0,0)$ with $a\ge 4$ is false by Lemma \ref{unbalanced}. To finish the proof, we use \thmref{main2} on $T(1,2,3;4;0,0,0)$ to show that $\sigma_4(X)$ fills the ambient space.
\qedd

\

\begin{prop0}  \label{carlini} If $X=\P {{1}} \times \P {{1}} \times \P {{n}} \times \P {{n}}$ then 
\begin{itemize}
\item[(i)] $X$ has a defective $2n+1$-secant variety.
\item[(ii)] The codimension of $\sigma_{2n+1}(X)$ is 2. 
\item[(iii)] $T(1,1,n,n;2n;0^4)$ and $T(1,1,n,n;2n+2;0^4)$ are true.
\end{itemize}
\end{prop0}

\proof\hskip -3pt . The proof of (i) follows an argument shown to us by Enrico Carlini (see also \cite{CGG4}). The proofs of (ii) and (iii) use the inductive method.

Proof of (i): Write $X$ as $(\P {{1}}\times \P {{n}})\times (\P {{1}}\times \P {{n}})$. Project the $2n+1$ points into each factor $(\P {{1}}\times \P {{n}})\subset \P {{2n+1}}$. Consider the hyperplanes $H_1,H_2$ in each $\P {{2n+1}}$ which pass through these projected points. Then the hyperplane defined by $H_1\otimes H_2$ contains the tangent space to $X$ at each of the $2n+1$ points. We can repeat this argument by switching the copies of $\P {{n}}$ to obtain a second pair of hyperplanes $H_1',H_2'$. Then the hyperplane defined by $H_1'\otimes H_2'$ also contains the tangent spaces to $X$ at each of the $2n+1$ points. Thus by Terracini's Lemma, the codimension of $\sigma_{2n+1}(X)$ is at least 2.

\

Proof of (ii): It is enough to show that $T(0,1,1,n+1,n+1;2n+3;2,0^4)$ is true. This is a superabundant case that reduces by \thmref{main2} to
$$T(0,1,1,n+1,n;2n+1;2,0^3,2)\hskip 6pt {\rm and} \hskip 6pt T(0,1,1,n+1,0;2;0,0^3,2n+1).$$ The second of these statements is true since no Segre variety has a defective 2-secant variety. Note that $(0,1,1,n+1,n;2n+1;2,0^3,2)$ is equiabundant so we use \corref{main5} to reduce $T(0,1,1,n+1,n;2n+1;2,0^3,2)$ to $T(1,1,n+1,n;2n+1;0^3,2)$ then use \thmref{main1} to reduce to
$$T(1,1,n,n;2n;0,0,1,1)\hskip 6pt {\rm and} \hskip 6pt T(1,1,0,n;1;0,0,2n,1).$$ The second of these statements is true. \thmref{main1} reduces
 $T(1,1,n,n;2n;0,0,1,1)$ to $$T(0,1,n,n;n;n,0,1,0)\hskip 6pt {\rm and} \hskip 6pt T(0,1,n,n;n;n,0,0,1).$$ These two statements are equivalent. \corref{main5} reduces $T(0,1,n,n;n;n,0,1,0)$ to $T(1,n,n;n;0,1,0)$ then we use \thmref{main1} to reduce to
$$n*T(1,n,0;1;0,0,n-1)\hskip 6pt {\rm and} \hskip 6pt T(1,n,0;0;0,1,n).$$ Both these statements are true so we are done.

\

Proof of (iii):  Since $T(1,1,n,n;2n;0^4)$ is subabundant, we use \thmref{main1} to reduce to $2*T(0,1,n,n;n;n,0,0,0)$. \corref{main5} reduces $T(0,1,n,n;n;n,0,0,0)$ to $T(1,n,n;n;0,0,0)$. Finally, we use \thmref{main1} to reduce $T(1,n,n;n;0,0,0)$ to $n*T(1,n,0;1;0,0,n-1)$ and $T(1,n,0;0,0,0,n)$. Both these statements are true. 

Since $T(1,1,n,n;2n+2;0^4)$ is superabundant, we use \thmref{main2} to reduce to $(n+1)*T(1,1,n,0;2;0,0,0,2n)$. This statement is true since no Segre variety has a defective 2-secant variety.
\qedd

\

\begin{remark} \propref{carlini} gives a complete description of the dimensions of the secant varieties to $X=\P {{1}} \times \P {{1}} \times \P {{n}} \times \P {{n}}$. In particular, $X$ has no defective $p$-secant varieties for $p\le 2n$ and $\sigma_{2n+2}(X)$ fills the ambient space, that is $\underline R(1,1,n,n)=2n+2$.
\end{remark}

It is interesting to compare the following proposition with \propref{carlini}

\begin{prop0}
For any positive integer $n$, $T(1,1,n,n+1;2(n+1),0^4)$ is perfect.
\end{prop0}

\proof The statement reduces to $T(1,1,0,n+1;2,0^2,2n,0)$ which is true because
 $T(1,1,n+1;2,0^3)$ is true and subabundant, and the $2n$ additional conditions are independent.\qedd

\

\begin{prop0}
\label{p2p3p3}
$\dim\sigma_5(\P {{2}}\times \P {{3}}\times\P {{3}})=43$.
\end{prop0}

\proof\hskip -3pt . We first show that $\P {{2}}\times \P {{3}}\times\P {{3}}$ has a defective 5-secant variety. In other words, we show that $\dim\sigma_5(\P {{2}}\times \P {{3}}\times\P {{3}})<44$. Given five general points in $X=\P {{2}}\times \P {{3}}\times\P {{3}}$ we want to construct a rational
normal curve of degree $8$, $C_8\subset X$, passing through the five points. We project the five points from $X$ onto each factor.
We get on $\P {{2}}$ a conic $C_2$ through five points $Q_1,\ldots Q_5$, and an isomorphism
$g\colon \P {{1}}\to C_2$ such that $g(0)=Q_1$, $g(1)=Q_2$, $g(\infty)=Q_3$, $g(x_1)=Q_4$, $g(x_2)=Q_5$
for some points $x_1, x_2\in\P {{1}}$.
In $\P {{3}}$ there is a two dimensional family of twisted cubics $C_{s,t}$ through the five projected points $P_1,\ldots P_5\in \P {{3}}$.
This means we have a family of maps $f_{s,t}\colon\P {{1}}\to C_{s,t}$ such that
$f_{s,t}(0)=P_1$, $f_{s,t}(1)=P_2$, $f_{s,t}(\infty)=P_3$. It is easy to see that
the preimage, $f_{s,t}^{-1}(P_4)$, is not constant when $s,t$ change. This fact can be verified by projecting from $P_5$ on a plane, where we get a pencil of conics through four points, and it is straightforward to check that the cross ratio of the four points is not constant in the pencil. Then we can choose $s$, $t$ such that
$f_{s,t}(x_1)=P_4$, $f_{s,t}(x_2)=P_5$. Repeating the same argument for the second copy of $\P {{3}}$ we get a morphism $f\colon\P {{1}}\rightarrow  \P {{2}}\times \P {{3}}\times\P {{3}}$ through
the five original points of degree $2+3+3=8$. This is the desired $C_8$ which spans
a space $\P {{8}}$. Hence each of the five tangent spaces
at the five original points meets this $\P {{8}}$ in a line and the span of the five tangent spaces has dimension $\le 8+5\cdot 7=43$. By Terracini's lemma this concludes the proof.

\

\begin{remark} $T(2,3,3;s;0^3)$ is true if $s\le 4$, moreover $\underline R(2,3,3)=6$.
\end{remark}

Now we show that $\dim\sigma_5(\P {{2}}\times \P {{3}}\times\P {{3}})=43$. It is enough to show that $T(0,2,3,3;5;4,0,0,0)$ is true. We use \thmref{main2} to reduce $T(0,2,3,3;5;4,0,0,0)$  to $T(0,2,3,0;1;2,0,0,4)$ and $T(0,2,3,2;4;2,0,0,1)$. The first of these statements is true. We use \thmref{main1} to reduce  $T(0,2,3,2;4;2,0,0,1)$ to $T(0,2,1,2;2;2,0,2,0)$ and $T(0,2,1,2;2;0,0,2,1)$. Both of these statements are true from \propref{listdef}. \qedd

\

\begin{thm0}
$\sigma_5(\P {{n_1}}\times\ldots\times\P {{n_k}})$  is non-defective with the following exceptions:
 
$(n_1, n_2, n_3)=(2,3,3)$  

$(n_1, n_2, n_3)=(1,2,a)$ with $a\ge 5$

$(n_1, n_2, n_3)=(1,3,a)$ with $a\ge 5$

$(n_1, n_2, n_3,n_4)=(1,1,2,2)$

\end{thm0}

\proof\hskip -3pt .  By \cite{CGG2}, $T(1,1,1,1,1;5;0,0,0,0,0)$ is known to be true.  Thus there are no exceptions with $k\ge 5$. To treat the case $k=4$ we prove that $T(1,1,1,4;5;0,0,0,0)$ is true. By \thmref{main1} we reduce $T(1,1,1,4;5;0,0,0,0)$ to $2*T(1,1,1,1;2;0,0,0,3)$ and $T(1,1,1,0;1;0,0,0,4)$. All these statements are known to be true. In the same manner we prove that $T(1,1,2,3;5;0,0,0,0)$ is true.
By \thmref{main2}, we show that $T(1,1,1,2;5;0,0,0,0)$ is true. By \propref{carlini},  $T(1,1,2,2;5;0,0,0,0)$ is false.

Now we treat the case $k=3$. Let us begin by proving that $T(2,2,4;5;0,0,0)$ is true. Indeed we reduce by \thmref{main1} to
$T(2,2,1;2;0,0,3)$ and $T(2,2,2;3;0,0,2)$ which are both true.
Similarly by \thmref{main2}, $T(2,2,3;5;0,0,0)$ is true. By \propref{p2p3p3}, $T(2,3,3;5;0,0,0)$ is false.
Hence if $n_1\ge 2$ the theorem is true and we may assume $n_1=1$.

Let us now prove that  $T(1,4,4;5;0,0,0)$ is true.
Note that the 7-tuple is equiabundant. We use \thmref{main1}  and reduce  to
$T(1,1,4;2;0,3,0)$ and $T(1,2,4;3;0,2,0)$
and again to
$T(1,1,2;1;0,2,1),\ T(1,1,1;1;0,1,1),\ T(1,2,2;2;0,1,1),\  T(1,2,1;1;0,1,2)$. All these statements are known to be true.
Hence if $n_2\ge 4$ the theorem is true and we may assume $n_2\le 3$. The cases $(n_1, n_2, n_3)=(1,2,a)$ with $a\ge 5$ and $(n_1, n_2, n_3)=(1,3,a)$ with $a\ge 5$ are defective by Lemma \ref{unbalanced}. To finish the proof, we note that $(1,3,4;5;0,0,0)$ is equiabundant and that $T(1,3,4;5;0,0,0)$ is true.
\qedd

\

\begin{thm0}
$\sigma_6(\P {{n_1}}\times\ldots\times\P {{n_k}})$  is non-defective with the following exceptions:
 
$(n_1, n_2, n_3)=(1,3,a)$  with $a\ge 6$

$(n_1, n_2, n_3)=(1,4,a)$ with $a\ge 6$

$(n_1, n_2, n_3)=(2,2,a)$ with $a\ge 6$

$(n_1, n_2, n_3,n_4)=(1,1,1,a)$ with $a\ge 6$
\end{thm0}

\proof\hskip -3pt .  The exceptions all follow from Lemma \ref{unbalanced}. To show there are no more exceptions, one needs to show that $T(\vec {\mathbf n},6,\vec {\mathbf 0})$ is true for the following values of $\vec {\mathbf n}$:

\

Subabundant cases: $(1^6), (1^4,2), (1,1,2,3), (1,2,2,2), (1,5,5), (3,3,3), (2,3,4)$

Superabundant cases: $(1^5), (1^3,5), (1,1,2,3), (1,4,5), (2,3,3), (2,2,5), (1,2,a)$

\

The subabundant cases can all be established using \thmref{main1} and \corref{main5}.
The case $(1,2,a)$ follows quickly from the case $(1,4,5)$ and \thmref{main2}. The other superabundant cases can all be established using \thmref{main2}.
\qedd

\

\section{Non-defectivity for many copies of $\P n$}
In this section we study Segre varieties of the form $X=\P {n} \times \dots \times \P {n}$.
We show that for most values of $s$, $\sigma_s(X)$ is non-defective. Before we prove the main theorem, we need a technical lemma.

\begin{lemma0}\label{tech} Let $\tilde{s_k}=\frac{(n+1)^k}{nk+1}$ and $s_k=\lfloor\tilde{s_k}\rfloor$. Let $\delta_k\equiv s_k\  mod\ (n+1)$  with $\delta_k\in\{0,\dots ,n\}$. Let $q=\frac{s_k-\delta_k}{n+1}$ and $\tilde {q}=\frac{\tilde{s_k}-\delta_k}{n+1}$.
\begin{itemize}
\item[(i)] If ($k=4$ and $n\ge 12$) or if ($k=5$ and $n\ge 4$) or if ($k=6,7$ or $8$ and $n\ge 2$) or if ($k\ge 9$ and $n\ge 1$) then $q+1\le s_{k-1}-\delta_{k-1}$.
\item[(ii)] $(q+1)(nk-n+1)+(s_k-\delta_k)-q+n\ge (n+1)^{k-1}$
\end{itemize}
\end{lemma0}

\proof\hskip -3pt . We write each proof as a sequence of implications. 

Proof of (i): The first statement follows from the fact that $q+1$ is an integer. \begin{eqnarray*}
& q+1\le s_{k-1}-\delta_{k-1}\\
  \iff & q+1\le \tilde{s}_{k-1}-\delta_{k-1}\\
\Longleftarrow\ & \tilde{q}+1\le \tilde{s}_{k-1}-\delta_{k-1}\\
\iff & \tilde{s}_k-\delta_k + (n+1) \le (n+1)(\tilde {s}_{k-1}-\delta_{k-1})\\
\iff & (n+1)\delta_{k-1} - \delta_k +(n+1)\le (n+1)\tilde {s}_{k-1}- \tilde {s}_k\\
\iff & (n+1)\delta_{k-1} - \delta_k +(n+1)\le(n+1)^k\left(\frac{1}{nk+1-n} - \frac{1}{nk+1}\right)\\
\iff & \frac{\delta_{k-1}}{n+1} - \frac{\delta_k}{(n+1)^2} + \frac{1}{n+1} \le(n+1)^{k-2}\left(\frac{1}{nk+1-n} - \frac{1}{nk+1}\right)\\
\iff & \frac{\delta_{k-1}}{n+1} - \frac{\delta_k}{(n+1)^2} + \frac{1}{n+1} \le(n+1)^{k-4}\left(\frac{n(n+1)^2}{(nk+1-n)(nk+1)}\right)\\
\end{eqnarray*}

Since $\delta_{k-1}\le n$, the last statement is implied by $$1\le (n+1)^{k-4}\left(\frac{n(n+1)^2}{(nk+1-n)(nk+1)}\right).$$ Now the conclusions of part (i) are easy.

\

Proof of (ii): \begin{eqnarray*}
& (q+1)(nk-n+1)+(s_k-\delta_k)-q+n\ge (n+1)^{k-1}\\
\iff & (n+1)[(q+1)(nk-n+1)+(s_k-\delta_k)-q+n]\ge (n+1)^k\\
\iff & (s_k-\delta_k)(nk+1) + (n+1)(nk+1)\ge (n+1)^k\\
\iff & (nk+1)(s_k-\delta_k+n+1)\ge (n+1)^k\\
\end{eqnarray*}
Now this last statement is implied by $$(nk+1)(\tilde{s}_k-\delta_k+n)\ge (n+1)^k$$ which is equivalent to $$(n-\delta_k)(nk+1)\ge 0.$$ Since $\delta_k\le n$, we are done.
\qedd

\

\begin{thm0} \label{tensorpower}
Let $X=(\P {{n}})^k$, $k\ge 3$. Let $s_k$ and $\delta_k$ be defined by
$$s_k=\left\lfloor\frac{(n+1)^k}{nk+1}\right\rfloor \hskip 8pt{\rm and}\hskip 8pt  \delta_k\equiv s_k\  mod\ (n+1) \hskip 8pt {\rm with}\hskip 8pt \delta_k\in\{0,\dots ,n\}.$$

(i) If $s\le s_k-\delta_k$ then $\sigma_s(X)$ has the expected dimension.

(ii) If $s\ge s_k-\delta_k+n+1$ then $\sigma_s(X)$ fills the ambient space.
\end{thm0}

\proof\hskip -3pt .   The proof is by induction on $k$.  

Proof of (i): Note that $(n^k;s_k-\delta_k; 0^k)$ is subabundant. We start from the fact that $(\P {n})^3$ is nondefective when $n\neq 2$ \cite{L} and the fact that  $(\P {2})^4$ is nondefective.
Suppose that $T(n^{k-1};s_{k-1}-\delta_{k-1}; 0^{k-1})$ is true with $k\ge 4$. We need to show that
$T(n^k;s_k-\delta_k; 0^k)$ is true. If $q=\frac{s_{k}-\delta_k}{n+1}$ then we use \thmref{main1} to reduce
$T(n^k;s_k-\delta_k; 0^k)$ to
$T(0,n^{k-1};q;(s_k-\delta_k-q), 0^{k-1})$. Since $(0,n^{k-1};q;(s_k-\delta_k-q), 0^{k-1})$ is subabundant, we can reduce $T(0,n^{k-1};q;(s_k-\delta_k-q), 0^{k-1})$ to $T(n^{k-1};q;0^{k-1})$. By the induction hypothesis  we have
$T(n^{k-1};s_{k-1}-\delta_{k-1}; 0^{k-1})$ is true. If we can show that $q\le s_{k-1}-\delta_{k-1}$ then we are done. By Lemma \ref{tech}, we have
$q\le s_{k-1}-\delta_{k-1}$ with a small number of possible exceptions. Using the exact inequality, the only true exceptions are
$(n,k)=(4,4)$ or $(n,k)=(7,4)$. With the aid of a computer, we can take care of these cases by showing that  $T(4^4;36;0^4)$ and $T(7^4;136;0^4)$ are true. 

\

Proof of (ii): Note that $(n^k;s_k-\delta_k+n+1; 0^k)$ is superabundant. We again start from the fact that $(\P {n})^3$ is nondefective when $n\neq 2$ and the fact that  $(\P {2})^4$ is nondefective.
Suppose that $T(n^{k-1};s_{k-1}-\delta_{k-1}+n+1; 0^{k-1})$ is true with $k\ge 4$. We need to show that
$T(n^k;s_k-\delta_k+n+1; 0^k)$ is true. We use \thmref{main2} to reduce $T(n^k;s_k-\delta_k+n+1; 0^k)$ to
$T(0,n^{k-1};q+1;s_k-\delta_k-q+n,0^{k-1})$. From the proof of the first part of this theorem, we know that $T(n^{k-1};q+1;0^{k-1})$ is true except for a small number of possible exceptions. Using the exact inequality, the only true exceptions are $(n,k)=(1,5),(1,6),(1,7),(2,4),(3,4),(3,5),(4,4),(7,4)$. With the aid of a computer, we can show that $(\P {n})^k$ is non-defective in all of these cases. Thus we know the dimension of the variety corresponding to $T(n^{k-1};q+1;0^{k-1})$. To establish that $T(0,n^{k-1};q+1;s_k-\delta_k-q+n,0^{k-1})$ is true, note that we have $s_k-\delta_k-q+n$ ``point conditions". Such conditions are always independent, they correspond to adding in $s_k-\delta_k-q+n$ general vectors before computing the span. We want to show that the partial secant variety corresponding to $T(0,n^{k-1};q+1;s_k-\delta_k-q+n,0^{k-1})$ fills the space. In other words, we need to show that $$(q+1)(nk-n+1)+(s_k-\delta_k)-q+n\ge (n+1)^{k-1}.$$ But this statement follows from Lemma \ref{tech}. 
\qedd

\

The following corollary applies in the cases considered in the prop 2.2. of \cite{LM}.

\begin{coro0}
$T((r-1)^k,r,0^k)$ is subabundant and true if $k\ge 3$, $\forall r\ge 1$. 
\end{coro0}

\begin{remark} The particular case when $n=1$ in Theorem \ref{tensorpower} appears as Theorem 2.3 in  \cite{CGG2}.
 It is worth emphasizing that Theorem \ref{tensorpower} states that  $X=(\P {n})^k$ has at most $n$ values of $s$
 for which $T(n^k;s;0^k)$ is not true.  In many cases the inequalities of the previous theorem can be improved
by looking at the arithmetic of the particular numbers involved. An example of this phenomenon can be seen in the following
 corollary and example.  See also \propref{p18} and \propref{p23} which show that in some cases
$X=(\P {n})^k$ has at most one value of $s$ for which $T(n^k;s;0^k)$ is not true.
\end{remark}

\begin{coro0}  \label{PP} If $X=(\P {n})^k$ is numerically perfect and $\delta_k=0$ then $X$ is perfect.
\end{coro0}

\begin{exa} We can apply Corollary \ref{PP} if and only if $\frac{(n+1)^{k-1}}{nk+1}$ is an integer. For instance, if $(n+1=p^h$ for some prime number $p)$ and
$(k=\frac{p^{th}-1}{p^h-1}$ for some $t\ge 2)$
then $X$ is perfect. This example appeared in \cite{CGG1} utilizing some ideas from coding theory.
\end{exa}

An easy consequence of \thmref{tensorpower} is
\begin{coro0}
$$\underline R(n^k) \sim  \frac{(n+1)^k}{nk+1}$$
when $n\to\infty$ or  $k\to\infty$.
\end{coro0}

Let's take a closer look at the case $X=(\P {3})^k$. 
Lickteig showed that $(\P 3)^3$ is non-defective. \corref{PP} shows 
that $(\P 3)^5$ is non-defective. 
According to \thmref{tensorpower} we have that $T(3^4;16;0^4)$  and 
$T(3^4,20;0^4)$ are true, in particular $\underline R(3^4)=20$.
We want to show that $T(3^4;18;0^4)$ is true.
This will show that the inductive technique often goes further than the statement of
\thmref{tensorpower}.
In order to study  $(\P 3)^4$ we will need the following lemma.

\begin{lemma0}
\label{auxiliar}
$T(1^4;2;1,1,0,0)$, $T(1^4;1;2,1,1,1)$ and $T(1^4;0;2^4)$ are true
\end{lemma0}

\proof\hskip -3pt .   Use Corollary \ref{main5} to reduce to $(\P 1)^3$.\qedd

\

\begin{prop0} \label{p18}
$T(3^4;18;0^4)$ is true, that is  $\sigma_{18}(\P 3)^4$ has the expected dimension.
\end{prop0}

\proof\hskip -3pt .   We use \thmref{main1} to reduce to two copies of $\P 1\times (\P 3)^3$, then
to four copies of $(\P 1)^2\times (\P 3)^2$, then to eight copies of $(\P 1)^3\times (\P 3)^1$,
then to sixteen copies of $(\P 1)^4$. In the end we need sixteen 5 tuples $(s,a_1,a_2,a_3,a_4)$
such that $T(1^4;s;a_1,a_2,a_3,a_4)$ is true and such that the vector sum of the sixteen 5 tuples
is $(18, 18, 18, 18, 18)$. Utilizing \lemref{auxiliar}, a solution is accomplished by the following eight vectors repeated twice.
$$(2,1,1,0,0)$$
$$(2,0,0,1,1)$$
$$(1,2,1,1,1)$$
$$(1,1,2,1,1)$$
$$(1,1,1,2,1)$$
$$(1,1,1,1,2)$$
$$(1,1,1,1,1)$$
$$(0,2,2,2,2)$$\qedd

\

We want to show that the inductive technique goes  further than the statement of
\thmref{tensorpower} also in the superabundant case.

Indeed we know that $\underline R(2^3)=5$ (defective case) and $\underline R(2^4)=9$
(\corref{PP}). According to \thmref{tensorpower} we have that
$T(2^5;21;0^5)$ is true and that $23\le \underline R(2^5)\le 24$.

We can show that

\begin{prop0} \label{p23}
$$\underline R(2^5)= 23$$
\end{prop0}

\proof\hskip -3pt . Since $T(2^5;22;0^5)$ is subabundant and true and since $(2^5;23;0^5)$ is superabundant,
 it is enough to show that $T(2^5;23;0^5)$ is true. We use \thmref{main2} to reduce to
$T(2^4,1;15;0^4,8)$ and $T(2^4,0;8;0^4,15)$. Since $(\P {2})^4$ is perfect, $T(2^4,0;8;0^4,15)$ is true.
 We use \thmref{main2} to reduce $T(2^4,1;15;0^4,8)$ to $T(2^3,0,1;5;0^3,10,2)$ and $2*T(2^3,0,1;5;0^3,10,3)$.
 Since $T(2^3,0,1;5;0^3,10,2)$ implies $T(2^3,0,1;5;0^3,10,3)$, it is enough to show that $T(2^3,0,1;5;0^3,10,2)$ is true.
 We use \corref{main6} to reduce $T(2^3,0,1;5;0^3,10,2)$ to $T(2^3,1;5;0^3,2)$. Now use \thmref{main1} to reduce
 to $T(2^3,0;2;0^3,5)$ and $T(2^3,0;3;0^3,4)$. Both of these statements are true from the classification
 of Segre varieties with defective 3-secant varieties. \qedd

\

We do not currently have a general theorem that shows that every tensor power of $\P {n}$ is non-defective. However, if $n$ is odd, we can prove that  for each tensor power of $\P {n}$, there exists a Segre product of a projective space with the tensor power which is not only non-defective but perfect.

\begin{thm0}\label{ktn}
If $n$ is odd then the Segre variety $\P {k}\times (\P {n})^{k+1}$ is perfect.
\end{thm0}

\proof\hskip -3pt .   First note that $T(k,n,\dots,n;(n+1)^k;0,\dots,0)$ is equiabundant. 

\noindent Since $n$ is odd, $$T(k,n,\dots,n;(n+1)^k;0,\dots,0)$$ reduces to (multiple copies of)
$$T(k,1,n,\dots,n;2(n+1)^{k-1};0,(n+1)^{k-1}(n-1),0,\dots,0)$$ and then to
$$T(k,1,1,n,\dots,n;4(n+1)^{k-2};0,2(n+1)^{k-2}(n-1),2(n+1)^{k-2}(n-1),0,\dots,0).$$
We continue in this manner until we reduce to
$$T(k,1,\dots,1,n;2^k;0,2^{k-1}(n-1),\dots,2^{k-1}(n-1),0).$$
Now we reduce to $\frac{n-1}{2}$ copies of 
$$T(k,1,\dots,1;0;0,2^k,\dots,2^k)$$ and one copy of
$$T(k,1,\dots,1;2^k;0,\dots,0).$$
Iterating \corref{main5} we reduce $T(k,1,\dots,1;0;0,2^k,\dots,2^k)$ to $T(k,1,1;0;0,2,2)$. In a similar manner we reduce $T(k,1,\dots,1;2^k;0,\dots,0)$ to $T(k,1,1;2;0,0,0)$. Both of these statements are true and we are done.
\qedd

\

\section{Closing Remarks and Open Questions}

\noindent CLASSIFICATION OF DEFECTIVE  $\sigma_s(X)$

 By Lemma \ref{unbalanced}, we know that unbalanced Segre varieties are defective. Using a Montecarlo technique combined with Terracini's Lemma (as in \cite{Mc}), we can show there are no balanced $t$-defective Segre varieties ($t\le 8$) other than the known  cases: $(2,2,2),(2,3,3),(2,4,4),(1,1,1,1),(1,1,2,2),(1,1,3,3)$. The cases $(2,2,2)$ and $(2,4,4)$ are in a family originally described by Strassen. The three cases $(1,1,1,1),(1,1,2,2)$ and $(1,1,3,3)$ are in the family covered by \propref{carlini}. The case $(2,3,3)$ seems to fall into its own family and is proven to be defective in \propref{p2p3p3}. Thus, all known cases of defective Segre varieties fall into one of the following four families: \{{\it unbalanced}, $(1,1,n,n), (2,3,3), (2,n,n)$ with $n$ even\}.
 
With the aid of a computer combined with a Montecarlo technique, we can show that every balanced, numerically perfect, 3 odd factor Segre Variety with $n_3\le 30$ is perfect. With the use of the inductive procedure combined with computer calculations, most of the balanced, numerically perfect cases with $n_3\le 100$ can be shown to be perfect.

\

\noindent MANY COPIES OF $\P {n}$

Arithmetical properties of $n$ and $k$ often allow \thmref{tensorpower} to be improved in special cases as we did in \propref{p18} and \propref{p23}. When $k\ge 3$, we strongly suspect there are only a finite number of defective Segre varieties of the form $(\P {n})^k$. We somewhat suspect that $(\P {2})^3$ and $(\P {1})^4$ are the only defective cases.

\

\noindent OPEN QUESTIONS

\begin{question} Let $X=\P {{n_1}}\times \P {{n_2}} \times \P {{n_3}}$. If $X$ is numerically perfect and balanced with $n_1,n_2,n_3$ odd then is $X$ perfect?
\end{question}

\begin{question} If $X=\P {{n_1}}\times \P {{n_2}} \times \P {{n_3}}$ is numerically perfect and balanced then is $X$ perfect?
\end{question}

\begin{question} Do all defective Segre varieties of the form $X=\P {{n_1}}\times \P {{n_2}} \times \P {{n_3}}$ fall into the following 3 classes:
\begin{itemize}
\item[1.] $X$ is unbalanced.
\item[2.] $X=\P {{2}}\times \P {{n}} \times \P {{n}}$ with $n$ even.
\item[3.] $X=\P {{2}}\times \P {{3}} \times \P {{3}}$.
\end{itemize}
\end{question}

\begin{question} Let $k\ge 3$. Other than $(\P {2})^3$ and $(\P {1})^4$, is every Segre variety of the form $(\P {n})^k$ nondefective?
\end{question}

\begin{question} Does  there exist a $T$ such that ${\bf P}^{\vec{\bf n}}$ is nondefective whenever $(n_1,\dots n_k)$ is balanced with $k>T$?
\end{question}

\begin{question} Do all defective Segre varieties fall into the following 4 classes:
\begin{itemize}
\item[1.] $X$ is unbalanced.
\item[2.] $X=\P {{2}}\times \P {{n}} \times \P {{n}}$ with $n$ even.
\item[3.] $X=\P {{2}}\times \P {{3}} \times \P {{3}}$.
\item[4.] $X=\P {{1}}\times \P {{1}} \times \P {{n}} \times \P {{n}}$.
\end{itemize}
\end{question}

\

\noindent {\it Acknowledgments}: We thank Luca Chiantini, Ciro Ciliberto and Tony Geramita for helpful and interesting discussions on secant varieties. We also thank Enrico Carlini for the proof of part (i) of Theorem \ref{carlini}.

\

\

\

\noindent {\it Hirotachi Abo}

\noindent Department of Mathematics, Colorado State University, Fort Collins, Colorado 80525

\noindent  abo@math.colostate.edu, http://www.math.colostate.edu/$\sim$abo 
\vskip 3pt
\noindent {\it New address: Department of Mathematics, University of Idaho, Moscow, Idaho 83844}

\

\noindent {\it Giorgio Ottaviani}

\noindent Dipartimento di Matematica ``Ulisse Dini", Universit\`a degli Studi di Firenze, Viale Morgagni 67/A, 50134 Firenze, Italy

\noindent ottavian@math.unifi.it, http://www.math.unifi.it/$\sim$ottavian

\

\noindent {\it Chris Peterson}

\noindent Department of Mathematics, Colorado State University, Fort Collins, Colorado 80525

\noindent peterson@math.colostate.edu, http://www.math.colostate.edu/$\sim$peterson

\end{document}